\documentclass{amsart}
\pdfoutput=1
  \usepackage[english]{babel}    

  \usepackage[latin1]{inputenc} 
  \usepackage{amsmath}
  \usepackage{amsthm} 
  \usepackage{amsfonts} 
  \usepackage{mathtools}
  \usepackage{relsize}   
  \usepackage{tikz}      
  \usepackage{newtxtext}
  \usepackage{newtxmath}
  \usepackage{tikz-cd}
  \usepackage{enumerate}
  \usepackage{pict2e}
  \usepackage[backend=biber]{biblatex}
  \PassOptionsToPackage{hyphens}{url}\usepackage{hyperref}
  \usepackage{url}
  \usepackage{csquotes}
  \usepackage{extarrows}

  \newcommand{\nospacepunct}[1]{\makebox[0pt][l]{\,#1}}

  \newcommand{\colim}{\operatornamewithlimits{colim}}

  \newcommand{\op}{^{\operatorname{op}}}
  \newcommand{\cts}{^{\operatorname{cts}}}
  \newcommand{\shape}[2]{\hat{\Pi}_{(\infty,1)}^{#1}({#2})}
  \newcommand{\hooklongrightarrow}{\lhook\joinrel\longrightarrow}
  \DeclareMathOperator{\id}{id}
  
  \DeclareMathOperator{\map}{map}
  
  \DeclareMathOperator{\Tw}{Tw}
  \DeclareMathOperator{\Fun}{Fun}
  
  \DeclareMathOperator{\Psh}{Psh}
  
  \DeclareMathOperator{\Cat}{Cat}
  \DeclareMathOperator{\Sh}{Sh}
  
  \DeclareMathOperator{\Ind}{Ind}

  \DeclareMathOperator{\Proj}{Proj}
  
  \DeclareMathOperator{\Gal}{Gal}
  \DeclareMathOperator{\pr}{pr}

  \DeclareMathOperator{\disc}{disc}
  \DeclareMathOperator{\Cocart}{Cocart}
  \DeclareMathOperator{\eff}{eff}
  \DeclareMathOperator{\Un}{Un}
  \DeclareMathOperator{\hyp}{hyp}
  \DeclareMathOperator{\Pro}{Pro}
  \DeclareMathOperator{\CoCart}{CoCart}
  
  \DeclareMathOperator{\proet}{pro\acute{e}t}
  \DeclareMathOperator{\et}{\acute{e}t}
  \DeclareMathOperator{\Lfib}{Lfib}
  \DeclareMathOperator{\Pyk}{Pyk}
  \DeclareMathOperator{\pyk}{pyk}
  \DeclareMathOperator{\IPyk}{\mathbf{Pyk}}
  \DeclareMathOperator{\Str}{\mathbf{Str}}
  \DeclareMathOperator{\Map}{Map}
  \DeclareMathOperator{\Cond}{Cond}
  \newcommand{\coh}{\operatorname{coh}}
  \newcommand{\dash}{\nobreakdash-}

  \makeatletter
  \newcommand{\adjunction}{\@ifstar\named@adjunction\normal@adjunction}
  \newcommand{\normal@adjunction}[4]{%
    #1\colon #2%
    \mathrel{\vcenter{%
      \offinterlineskip\m@th
      \ialign{%
        \hfil$##$\hfil\cr
        \longrightharpoonup\cr
        \noalign{\kern-.3ex}
        \smallbot\cr
        \longleftharpoondown\cr
      }%
    }}%
    #3 \noloc #4%
  }
  \newcommand{\named@adjunction}[4]{%
    #2%
    \mathrel{\vcenter{%
      \offinterlineskip\m@th
      \ialign{%
        \hfil$##$\hfil\cr
        \scriptstyle#1\cr
        \noalign{\kern.1ex}
        \longrightharpoonup\cr
        \noalign{\kern-.3ex}
        \smallbot\cr
        \longleftharpoondown\cr
        \scriptstyle#4\cr
      }%
    }}%
    #3%
  }
  \newcommand{\longrightharpoonup}{\relbar\joinrel\rightharpoonup}
  \newcommand{\longleftharpoondown}{\leftharpoondown\joinrel\relbar}
  \newcommand\noloc{%
    \nobreak
    \mspace{6mu plus 1mu}
    {:}
    \nonscript\mkern-\thinmuskip
    \mathpunct{}
    \mspace{2mu}
  }
  \newcommand{\smallbot}{%
    \begingroup\setlength\unitlength{.15em}%
    \begin{picture}(1,1)
    \roundcap
    \polyline(0,0)(1,0)
    \polyline(0.5,0)(0.5,1)
    \end{picture}%
    \endgroup
  }
\makeatother

  \newtheorem{satz}{Theorem}[section]
  \newtheorem{lemma}[satz]{Lemma}
  \newtheorem{korollar}[satz]{Corollary}
  \newtheorem{proposition}[satz]{Proposition}
  \theoremstyle{definition}
  \newtheorem{definition}[satz]{Definition}
  \newtheorem{recollection}[satz]{Recollection}

  \newtheorem{construction}[satz]{Construction}
  
  \newtheorem{bemerkung}[satz]{Remark}
  \newtheorem{nix}[satz]{}
  \newenvironment{beweis}
    {\begin{proof}[Proof:]}
    {\end{proof}}

   \newenvironment{beweisv}[1]
   {\begin{proof}[Proof of {#1}:]}
   {\end{proof}}
  \newtheorem{beispiel}[satz]{Example}

  \newenvironment{sloppypar*}
  {\sloppy\ignorespaces}
  {\par}
  
  \numberwithin{equation}{satz}
  
  \tikzcdset{
  	cells={font=\everymath\expandafter{\the\everymath\displaystyle}},
  }

\begin{document}

\title{The pro-\'etale topos as a category of pyknotic presheaves}
\author{Sebastian Wolf}

\begin{abstract}
	In this paper we will show that for a quasicompact quasiseparated scheme $X$ the hypercomplete pro-\'etale $\infty$-topos, as introduced by Bhatt and Scholze, is equivalent to the $ \infty $-category of continuous representations of the Galois category $\Gal(X)$ of $X$ with values in the~$\infty$-category of pyknotic spaces.
	In particular this proves that internally to pyknotic spaces, the hypercomplete pro-\'etale $\infty$-topos of $X$ is an $ \infty $-category of presheaves.
\end{abstract}

\maketitle

\setcounter{tocdepth}{1}
\tableofcontents

\setcounter{section}{0}

\section{Introduction}

\subsection*{Main results}

For a quasicompact quasiseparated scheme $X$, the Galois category $\Gal(X)$ of $X$, as introduced in \cite{barwick2018exodromy}, is defined to be the $ \infty $-category of points of the \'etale $ \infty $-topos $X_{\et}$ of $X$.
This $ \infty $-category turns out to be a $ 1 $-category and admits a rather concrete description: 
The objects are given by geometric points of the scheme $X$ and the morphisms by specializations of such.
In \cite{barwick2018exodromy}, the authors refined $ \Gal(X) $ to a pro-$ \infty $-category whose underlying $ \infty $-category recovers the $ \infty $-category of points of $ X_{\et} $.
We can therefore make use of the ideas of \emph{condensed} or \emph{pyknotic mathematics} introduced in \cite{barwick2019pyknotic},  \cite{clausen2019condensed}, \cite{lurie2018ultracategories} and consider $ \Gal(X) $ as a pyknotic $ \infty $-category.
Recall that a pyknotic $ \infty $-category is defined to be a hypersheaf of $ \infty $-categories on the site of profinite sets and we can consider $ \Gal(X) $ as such via the assertion
\[
	K \mapsto \Fun(K,\Gal(X))
\]
for a pro-finite set $ K $.
Here we denote by $\Fun(K,\Gal(X))$ the $ \infty $-category of functors of pro-$ \infty $-categories from $ K $ to $ \Gal(X) $.

Given two pyknotic $ \infty $-categories $ \mathcal{C} $ and $ \mathcal{D} $, we can from the $ \infty $-category $ \Fun \cts (\mathcal{C},\mathcal{D}) $ of \emph{continuous} funtors between them (see Definition~\ref{def:FunctorCat}).
Let $ \mathcal{S}_{\pi} $ denote the $ \infty $-category of $ \pi $-finite spaces.
We can regard $ \mathcal{S}_{\pi} $ as a \emph{discrete} pyknotic $ \infty $-category $ \mathcal{S}_{\pi}^{\disc} $, by considering the constant sheaf on profinite sets with value $ \mathcal{S}_\pi $. 
In this language the \textit{exodromy equivalence} (see \cite[Corollary~0.5.2, Corollary~13.6.3]{barwick2018exodromy}) says that there is a natural equivalence
\[
	\Fun \cts ( \Gal(X),\mathcal{S}_\pi^{\disc}) \simeq X_{\et}^{\operatorname{constr}}
\]
between the $ \infty $-category of continuous functors from $ \Gal(X) $ to~$ \mathcal{S}_{\pi}^{\disc} $ and the $ \infty $-category of constructible \'etale sheaves on $ X $.

As explained in \cite{barwick2018exodromy}, the above equivalence holds in greater generality in the setting of so called \emph{spectral $ \infty $-topoi} (see \cite[Theorem~0.4.7, Definition 9.2.1] {barwick2018exodromy}).
These are $ \infty $-topoi $ \mathcal{X} $, equipped with a geometric morphism $ f_* \colon \mathcal{X} \to \Sh(S)$, where $ S $ is a spectral topological space, such that $ f_* $ exhibits the $ \infty $-category $ \operatorname{Pt}(X)$ as an $ S $-stratified space (see \cite[Proposition 9.2.5]{barwick2018exodromy} for the precise statement).
In this case $ \mathcal{X}$ is called a \emph{spectral $ S $-stratified $ \infty $-topos.}
The main example to keep in mind is the canonical geometric morphism $ X_{\et} \to X_{\operatorname{Zar}} $ for a qcqs scheme $ X $, from $ X_{\et} $ to the $ \infty $-topos $ X_{\operatorname{Zar}} $  of Zariski-sheaves on $ X $.
To any spectral $ S $-stratified $ \infty $-topos $ \mathcal{X} $, one can associate a pro-$ \infty $-category $ \shape S {\mathcal{X}} $, called the \emph{profinite $ S $-stratified shape of} $ \mathcal{X} $, whose underlying $ \infty $-category recovers the $ \infty $-category of points of $ \mathcal{X} $.
In this setting Barwick-Glasman-Haine show that there is an equvialence
\[
	\Fun \cts \big(\shape S {\mathcal{X}}, \mathcal{S}_\pi^{\disc}\big) \simeq \mathcal{X}^{S\operatorname{-constr}} 
\]
	between the $\infty$-category of continuous representations of the profinite $ S $-stratified shape $\shape S {\mathcal{X}}$ of $\mathcal{X}$ and the $ \infty $-category $ \mathcal{X}^{S\operatorname{-constr}} $ of $S$-constructible objects in $\mathcal{X}$ \cite[Theorem 0.4.7 and Definition 9.4.7]{barwick2018exodromy}.
	In the case $ \mathcal{X} = X_{\et}$, this recovers the previous equivalence.
	
	The goal of this paper is to extend the results above to a larger class of sheaves. 
	For this we recall from \cite[Construction 4.1.4]{barwick2019pyknotic} that we can equip the $ \infty $-category of \emph{pyknotic spaces} $ \Pyk(\mathcal{S}) $ with the structure of a pyknotic $ \infty $-category via the assertion
	\[
		K \mapsto \Pyk(\mathcal{S})_{/K}.
	\]
	We denote the resulting pyknotic $ \infty $-category by $ \IPyk(\mathcal{S}) $.
	The main theorem of this paper identifies the $ \infty $-category of continuous representations of $ \shape S {\mathcal{X}} $ with values in $ \IPyk(\mathcal{S}) $ as follows:

\begin{satz}\label{Theorem: Cts Exodromy}
	Let $ S $ be a spectral topological space and $\mathcal{X}$ an $S$-stratified spectral $\infty$-topos. 
	Then the exodromy equivalence induces an equivalence of $\infty$-topoi
	\[
		\mathcal{X}^{\pyk} \xlongrightarrow{\simeq} \Fun^{\operatorname{cts}}\big(\shape S {\mathcal{X}}, \mathbf{Pyk}(\mathcal{S})\big)
	\]
	between the pyknotification $\mathcal{X}^{\pyk}$ of $\mathcal{X}$ and the $ \infty $-category of continuous functors from the profinite stratified shape $\shape S {\mathcal{X}}$ to the $ \infty $-category of pyknotic spaces.
\end{satz}
Here the pyknotification \footnote{The pyknotification was introduced in \cite[Construction 3.3.2]{barwick2019pyknotic} under the name \emph{solidification}.} is defined as the $ \infty$-category of hypersheaves on the category of pro-$ S $-constructible sheaves $ \Pro(\mathcal{X}^{S\operatorname{-constr}}) $ where the topology is generated by finite families that are jointly effective epimorphisms (see Definition \ref{def:Pyknotification}).
The $ 1 $-categorical analogue of this construction was already considered in \cite[\S 7.1]{lurie2018ultracategories}. It follows from Example 7.1.7 of loc. cit. that for a qcqs scheme $ X $, the pyknotification $X_{\et}^{\pyk}$ is equivalent to the $ \infty $-topos of pro-\'etale hypersheaves $ X_{\proet}^{\hyp} $ on $ X $ as defined by Bhatt-Scholze in \cite{bhatt2013pro}.
This means that the process of passing from a spectral $ \infty $-topos to its pyknotification can be viewed as a generalization of the passage from the \'etale to the pro-\'etale $ \infty $-topos for a scheme $ X $.
Combining this alternative description of $ X_{\proet}^{\hyp} $ with Theorem~\ref{Theorem: Cts Exodromy} we arrive at the following Corollary:

\begin{korollar}\label{cor:X proet= Psh(Gal(X)op)}
	Let $X$ be a qcqs scheme. Then the exodromy equivalence induces an equivalence of $\infty$-topoi
	\[
		X_{\proet}^{\hyp} \xlongrightarrow{\simeq} \Fun \cts (\Gal(X), \IPyk (\mathcal{S})).
	\]
\end{korollar}
In particular, the above theorem says that from a pyknotic point of view, the hypercomplete pro-\'etale $ \infty $-topos is naturally equivalent to the $ \infty $-category of internal copresheaves on the pyknotic $ \infty $-category $\Gal(X)$.

This result has several interesting implications, some of which we try to explain in this article. 
First of all, many of the convenient properties of the pro-\'etale topos are now purely formal consequences of Theorem~\ref{Theorem: Cts Exodromy}.
For example the fact that the pro-\'etale topos is locally weakly contractible \cite[Proposition 4.2.8]{bhatt2013pro} can easily be deduced using our main result (see section 5).
Also it is a direct corollary of Theorem~\ref{Theorem: Cts Exodromy} that, contrary to the \'etale setting, for \emph{any} morphism $f \colon X \rightarrow Y$ of schemes the functor
\[
	f^* \colon Y_{\proet}^{\hyp} \rightarrow X_{\proet}^{\hyp}
\]
has the surprising property that it commutes with \emph{all} small limits (Corollary~\ref{cor: pullback pre limits}).
This generalizes \cite[Corollary 6.1.5]{bhatt2013pro} where this was shown for closed immersions only.
Our main theorem furthermore gives a conceptual explanation why the pro-\'etale topos has so many convenient categorical properties.
It is, in the pyknotic sense, just an $ \infty $-category of presheaves and thus inherits many nice properties of the much simpler $\infty$-topos $\Pyk(\mathcal{S})$ of pyknotic spaces. 
Furthermore this opens the way to understanding the pro-\'etale topos using \emph{pyknotic higher category theory} (i.e. higher category theory internal to the $\infty$-topos $\Pyk(\mathcal{S})$).
We will develop the basics of higher category theory internal to an $\infty$-topos in ongoing joint work with Louis Martini \cite{martini2022cocartesian,martini2021yoneda,martini2021limits,martini2022presentable}.

The language of \textit{pyknotic} or \textit{condensed} mathematics which we use in this article was developed by Barwick and Haine in \cite{barwick2019pyknotic} and independently at the same time by Clausen and Scholze in \cite{clausen2019condensed}.
The only real difference between the two approaches is the way in which set-theoretic issues are dealt with. 
Since we will directly build on the results in \cite{barwick2018exodromy} which are all formulated in the framework of \cite{barwick2019pyknotic}, we will do the same in this paper.
The 1-categorical versions of many ideas and constructions in this paper already appear in \cite[\S 6 and 7]{lurie2018ultracategories}.
Most notably the pyknotic category $ \Gal(X) $ is already introduced there under the name $ \operatorname{Stone}_{X_{\et}^{\operatorname{constr}}} $.
This construction is then used to prove an analogue of the exodromy equivalence in this setting. See \cite[Theorem 2.2.2]{lurie2018ultracategories} and \cite[\S 4.3]{barwick2019pyknotic}.

\subsection*{Technical Outline}
In the second section we will study categories of pro-objects and in particular effective epimorphisms in them.
This will allow us to construct the \textit{effective epimorphism topology} on $\Pro(\mathcal{X}_0)$ for a bounded $\infty$-pretopos $\mathcal{X}_0$ and thus define the \textit{pyknotification} of a bounded coherent $\infty$-topos.
It is quite elementary and may be mostly skipped by experts.

The main technical result heavily used in the proof of Theorem~\ref{Theorem: Cts Exodromy} is an alternative description of the $ \infty $-category of continuous functors
\[
	\Fun^{\operatorname{cts}}(\mathcal{C}, \mathbf{Pyk}(\mathcal{S}))
\]
for a pyknotic $\infty$-category $\mathcal{C}$. 
Namely we will see that it is equivalent to a suitably defined $ \infty $-category of \textit{continuous left fibrations} over $\mathcal{C}$ (see Theorem~\ref{theorem: continuous left fib are cts copresheaves}).
Providing this description and exploring a few consequences will be subject of the third section. 

In the fourth section we will proof Theorem~\ref{Theorem: Cts Exodromy}. 
The rough idea will be that the Exodromy equivalence together with the embedding $\mathcal{S}_\pi \hookrightarrow \Pyk(\mathcal{S})$ will provide an embedding 
	\[
		\mathcal{X}^{S\operatorname{-constr}} \simeq \Fun \cts\big(\shape S {\mathcal{X}},\mathcal{S}_\pi\big) \hooklongrightarrow \Fun \cts\big(\shape S {\mathcal{X}}, \IPyk (\mathcal{S})\big),
	\]
which will extend to a fully faithful (see Proposition~\ref{Prop: Pro(-) embeds fully faithful}) embedding
	\[
		\Pro(\mathcal{X}^{S\operatorname{-constr}}) \hooklongrightarrow \Fun \cts\big(\shape S {\mathcal{X}}, \IPyk (\mathcal{S})\big).
	\]
Since the pyknotification is by definition given by taking sheaves with respect to the effective epimorphism topology on $\Pro(\mathcal{X}^{S\operatorname{-constr}})$,
we would now like to apply Theorem~\ref{nices Sag theorem}, which would prove that $\Pro(\mathcal{X}^{S\operatorname{-constr}})$ equipped with the effective epimorphism topology is also a generating site for the right hand side above.
For this we have to show that the assumptions of Theorem~\ref{nices Sag theorem} are satisfied which will occupy the vast majority of this section.

\subsection*{Set theoretic conventions}

We will follow the set-theoretic conventions of \cite{barwick2019pyknotic}. 
In particular we fix a \textit{tiny} and a \textit{small} universe, respectively  determined by two strongly inaccessible cardinals
\[
	\delta_0 < \delta_1.
\]
We recall that, for a strongly inaccessible cardinal $\delta$, a $\delta$-$\infty$-topos is a left exact accessible localization of $\Fun(\mathcal{C}, \mathcal{S}_\delta)$ for a $\delta$-small $\infty$-category $\mathcal{C}$. 
We will simply say $\infty$-topos instead of~$\delta_1$\dash$\infty$\dash topos.

By convention, a \textit{spectral $S$-stratified $\infty$-topos} $\mathcal{X}$ will by always be a~$\delta_0$\dash$\infty$\dash topos.
In particular, the $\infty$-category $\mathcal{X}^{S\operatorname{-constr}}$ of $S$-constructible sheaves will always be tiny.

Recall that the $\infty$-category of \textit{pyknotic spaces} $ \Pyk(\mathcal{S}) $ is defined to be the $\infty$-topos of hypersheaves with values in $ \mathcal{S}_{\delta_1} $ on the site of tiny profinite sets equipped with finite jointly epimorphic families as coverings.
Any profinite set can be covered by an \emph{extremally disconnected set} $ K $ and these have the property that any epimorphism $ T \twoheadrightarrow K $ of profinite sets has a section.
Writing $ \Proj $ for the category of tiny extremally disconnected sets we therefore have an equivalence
\[
	\Pyk(\mathcal{S}) = \Fun^\times\big(\Proj \op, \mathcal{S}_{\delta_1}\big).
\]
It follows that $ \Pyk(\mathcal{S}) $ is hypercomplete and even postnikov-complete.
See \cite[\S 2.2]{barwick2019pyknotic} for more background on pyknotic spaces.

\subsection*{Acknowledgements}
I am very grateful to my advisor Denis-Charles Cisinski for many helpful discussions and for reading a preliminary draft.
Of course this work is massively based on the recent developments in \cite{barwick2018exodromy} due to Clark Barwick, Peter Haine and Saul Glasman and this paper may be viewed as a small addition to their work.
I would like to thank Clark Barwick especially, who already knew about Theorem~\ref{Theorem: Cts Exodromy} and encouraged me to work out a proof.
I would also like to thank the anonymous referee whose many suggestions greatly helped to improve the readability of this paper. 
The author was supported by the SFB 1085 `Higher Invariants' in Regensburg, funded by the DFG.

\section{Pro-objects and Pyknotification}

\begin{sloppypar}
The majority of the material presented in this section has been worked out in the 1\dash categorical case by  Lurie in \cite[\S 6.1]{lurie2018ultracategories} and most of the arguments in this section are just very straight-forward adaptions of the ones presented there.
\end{sloppypar}

\begin{recollection}
	Let $\mathcal{C}$ be a tiny $\infty$-category.
	We define the $ \infty $-category $ \Pro(\mathcal{C}) $ of \textit{pro-objects in $\mathcal{C}$} to be the essentially unique $ \infty $-category equipped with a functor $ i \colon \mathcal{C} \to \Pro(\mathcal{C}) $ such that for any $ \infty $-category with tiny cofiltered limits $ \mathcal{D} $, precomposition with $ i $ induces an equivalence
	\[
		\Fun'(\mathcal{C},\mathcal{D}) \xrightarrow{i^*} \Fun(\mathcal{C},\mathcal{D}).
	\]
	Here $ \Fun'(\mathcal{C},\mathcal{D}) $ denotes the $ \infty $-category of tiny cofiltered limit preserving functors.
	The $ \infty $-category $ \Pro(\mathcal{C}) $ exists by the dual of \cite[Proposition 5.3.6.2]{lurie2009higher} and is equivalent to $ \Ind(\mathcal{C} \op) \op $.
	In particular $ \Pro(\mathcal{C}) $ is small and locally tiny but in general not tiny.
\end{recollection}

\begin{recollection}
	Recall from \cite[Definition A.6.1.1]{lurie2016spectral} that an $ \infty $-category $ \mathcal{X}_0 $ is called an \emph{$ \infty $-pretopos} if
	\begin{enumerate}
		\item The $ \infty $-category $ \mathcal{X}_{0} $ has finite limits.
		\item Finite coproducts exist in $ \mathcal{X}_0 $ and are universal and disjoint.
		\item Groupoid objects in $ \mathcal{X}_0 $ are effective and their geometric realizations are universal.
	\end{enumerate}
	An $ \infty $-pretopos $ \mathcal{X}_0 $ is called \emph{bounded} if $ \mathcal{X}_0 $ is small and every object in $ \mathcal{X}_0 $ is $ n $-truncated for some integer $ n $. 
	If $ \mathcal{X} $ is a \emph{coherent} $ \infty $-topos in the sense of \cite[Definition A.2.0.12]{lurie2016spectral}, the full subcategory of $ \mathcal{X}^{\coh}_{<\infty} $ spanned by the truncated and coherent objects is a bounded $ \infty $-pretopos \cite[Example 7.4.4]{lurie2016spectral}.
	In fact any bounded $ \infty $-pretopos arises this way by \cite[Theorem 7.5.3]{lurie2016spectral}.
	For more background on $ \infty $-pretopoi and coherent $ \infty $-topoi, the reader may consult \cite[Appendix A]{lurie2016spectral} or \cite[\S 3]{barwick2018exodromy}
\end{recollection}

\begin{beispiel}
	If $ X $ is a qcqs scheme, the $ \infty $-topos $ X_{\et} $ of \'etale sheaves on $ X $ is coherent.
	In this case the $ \infty $-category $ X_{\et,<\infty}^{\coh} $ of truncated coherent objects is equivalent to the $ \infty $-category $ X_{\et}^{\operatorname{constr}} $ of constructible étale sheaves on $ X $, which is therefore an $ \infty $-pretopos.
	In fact we more generally have an equivalence
	\[
		\mathcal{X}^{\coh}_{<\infty} \simeq \mathcal{X}^{S\operatorname{-constr}}
	\]
	for any spectral $ S $-stratified $ \infty $-topos $ \mathcal{X} \to \Sh(S)$ \cite[Corollary 9.5.5]{barwick2018exodromy}.
\end{beispiel}

If $ \mathcal{X}_0 $ is an $ \infty $-pretopos, the $ \infty $-category $ \Pro(\mathcal{X}_0) $ is in general not an $ \infty $-pretopos.
The goal of this section is to show that even though this is case, the notion of effective epimorphism still yields a reasonable Grothendieck-topology on $ \Pro(\mathcal{X}_0) $. 

\begin{nix}
	Let $\mathcal{C}$ be an $\infty$-category with finite limits and geometric realizations of groupoid objects.
	Recall that a morphism $f \colon X \rightarrow Y$ in $\mathcal{C}$ is called an \textit{effective epimorphism}, if the canonical morphism
	\[
	\lvert \check{C}(f)_\bullet \rvert \longrightarrow Y
	\]
	is an equivalence, where $\check{C}(f)_\bullet$ denotes the \v{C}ech-nerve of $f$.
\end{nix}

\begin{lemma}\label{Lemma: limits of eff epis}
	Let ${\mathcal{X}_0}$ be a tiny bounded $\infty$-pretopos. 
	Let $f_\bullet \colon I \rightarrow \Fun(\Delta^1, \Pro({\mathcal{X}_0}))$ be a tiny cofiltered diagram of effective epimorphisms. 
	Then $\lim_i f_i$, considered as a morphism in~$\Pro({\mathcal{X}_0})$, is an effective epimorphism.
\end{lemma}
\begin{beweis}
	This is a straight-forward adaption of the argument given in \cite[Prop. E.5.5.3]{lurie2016spectral}:
	Let us denote the source and target of $f_\bullet$ by $X_\bullet$ and $Y_\bullet$, respectively.
	Let us write $U_\bullet$ for the \v{C}ech-nerve of $f= \lim_i f_i$.
	We would like to show that, for every $C \in \Pro({\mathcal{X}_0})$, the induced morphism
	\[
		\map_{\Pro({\mathcal{X}_0})}(\lim_i X_i, C) \longrightarrow \lim_{n \in \mathbf{\Delta}} \map_{\Pro({\mathcal{X}_0})}(U_n,C)
	\]
	is an equivalence.
	We observe that we may immediately assume that $C \in {\mathcal{X}_0} \subseteq \Pro({\mathcal{X}_0})$.
	We now write $U_{\bullet, i}$ for the \v{C}ech-nerve of $f_i$.
	Since we assumed that $C$ is cocompact, the above map may be identified with the composite
	\begin{align*}
		\colim_i \map_{\Pro(\mathcal{X}_0)}( X_i, C) &\longrightarrow \colim_i \lim_{n \in \mathbf{\Delta}}\map_{\Pro(\mathcal{X}_0)}(U_{i, n},C) \\ & \xlongrightarrow{\alpha} \lim_{n \in \mathbf{\Delta}} \colim_i \map_{\Pro(\mathcal{X}_0)}(U_{i,n},C).
	\end{align*}  
	Since every $f_i$ was assumed to be an effective epimorphism, the first map is an effective epimorphism. 
	Thus it suffices to see that $\alpha$ is an equivalence.
	Now, since ${\mathcal{X}_0}$ is bounded, there is an $n \in \mathbb{N}$ such that $C$ is $n$-truncated. 
	We may thus replace the $U_{i,n}$ with their~$n$\dash truncations~$\tau_{\leq n}(U_{i,n})$.
	It follows from \cite[Proposition A.1]{hesselholt2022dirac} that in the commutative diagram
	\[
		\begin{tikzcd}
{\colim_i \lim_{n \in \mathbf{\Delta}}\map_{\Pro(\mathcal{X}_0)}(U_{i, n},C)} \arrow[d] \arrow[r,"\alpha"] & {\lim_{n \in \mathbf{\Delta}} \colim_i \map_{\Pro(\mathcal{X}_0)}(U_{i,n},C)} \arrow[d] \\
{\colim_i \lim_{n \in \mathbf{\Delta}^{\leq n}}\map_{\Pro(\mathcal{X}_0)}(U_{i, n},C)} \arrow[r]  & {\lim_{n \in \mathbf{\Delta}^{\leq n}} \colim_i \map_{\Pro(\mathcal{X}_0)}(U_{i,n},C)} 
\end{tikzcd}
	\]
	the horizontal arrows are equivalences and the bottom vertical arrow is an equivalence as well since taking limits over $\Delta^{\leq n}$ commutes with filtered colimits.
	This completes the proof.
\end{beweis}

Let us quickly recall the following from \cite[Proposition A.6.2.1]{lurie2016spectral}:

\begin{proposition}
	Let ${\mathcal{X}_0}$ be an $\infty$-pretopos. 
	Then ${\mathcal{X}_0}$ admits a factorization system $(S_L,S_R)$ (in the sense of \cite[\S 5.2.8]{lurie2009higher}), where $S_L$ is the collection of effective epimorphisms and~$S_R$ the collection of $(-1)$-truncated morphisms.
\end{proposition}

\begin{proposition}\label{Proposition: Eff epis sind super in Pro(-)}
	Let ${\mathcal{X}_0}$ be a tiny bounded $\infty$-pretopos. Then the following hold:
	\begin{enumerate}
		\item The collection of effective epimorphisms and $(-1)$-truncated morphisms form a factorization system on $\Pro({\mathcal{X}_0})$. 
		\item A morphism in $\Pro({\mathcal{X}_0})$ is an effective epimorphism if and only if it can be written as a tiny inverse limit of effective epimorphisms in ${\mathcal{X}_0}$.
		\item Effective epimorphisms are stable under pullback in $\Pro({\mathcal{X}_0})$.
	\end{enumerate}
\end{proposition}
\begin{beweis}
	It is clear that effective epimorphisms and $(-1)$-truncated morphisms are stable under retracts.
	Furthermore we observe that effective epimorphisms are left orthogonal to~$(-1)$\dash truncated morphisms in any $\infty$-category with finite limits and geometric realizations.
	So let $f \colon X \rightarrow Z$ be a morphism in $\Pro({\mathcal{X}_0})$.
	By (the dual of) \cite[Proposition 4.2.2]{lurie2016spectral} we can write $ f \simeq \lim_i h_i \circ \lim_j g_j$, where $ h_i $ is a $ (-1) $-truncated morphism in $ \mathcal{X}_0 $ and $ g_j $ is an effective epimorphism in $ \mathcal{X}_0 $.
	By Lemma~\ref{Lemma: limits of eff epis} it follows that $\lim_i g_i$ is an effective epimorphism.
	Furthermore the diagonal $ \lim _i Y_i \rightarrow \lim _i Y_i \times_{\lim _i Z_i} \lim _i Y_i$ may be identified with the limit of the diagonals $Y_i \rightarrow Y_i \times_{Z_i} Y_i$ and is thus an equivalence.
	This proves i).
	
	One direction of ii) is simply Lemma~\ref{Lemma: limits of eff epis} combined with the observation that the inclusion $\mathcal{X}_0 \hookrightarrow \Pro(\mathcal{X}_0)$ preserves effective epimorphisms. 
	For the other direction assume that~$f \colon X \rightarrow Z$ is an effective epimorphism.
	We now pick a factorization $f \simeq (\lim_i  h_i) \circ (\lim_i g_i)$ as above. 
	By i) it follows from \cite[Proposition 5.2.8.6]{lurie2009higher} that $\lim_i h_i$ is both an effective epimorphism and $(-1)$-truncated.
	Thus it is left orthogonal to itself, so it is an equivalence, which completes the proof of ii).
	
	For iii) we may use \cite[Proposition 5.3.5.15]{lurie2009higher} again to assume that we are given a cofiltered diagram $I \rightarrow \Fun(\Lambda^2_2,{\mathcal{X}_0})$ depicted as 
	\[
		\begin{tikzcd}
                    & X_\bullet \arrow[d,"f_\bullet"] \\
T_\bullet \arrow[r,"\gamma_\bullet"] & Z_\bullet         
\end{tikzcd}
	\]
	such that the induced map $\lim_i X_i \rightarrow \lim_i Z_i$ is an effective epimorphism. We have to show that the induced map 
	\[
		\lim_i T_i \times_{\lim_i Z_i} \lim_i X_i \longrightarrow \lim_i T_i
	\]
	is an effective epimorphism.
	Again we get a functorial factorization
	\[
		X_i \xlongrightarrow{g_i} Y_i \xlongrightarrow{h_i} Z_i
	\]
	where $g_i$ is an effective epimorphism and $h_i$ is $(-1)$-truncated for all $i$.
	We now consider the diagram
	\[
		\begin{tikzcd}
\lim_i (T_i \times_{Z_i} X_i) \arrow[r] \arrow[d, "g'"]             & \lim_i X_i \arrow[d, "g= \lim_i g_i"] \\
\lim_i (T_i \times_{ Z_i} Y_i) \arrow[d, "h'"] \arrow[r] & \lim_i Y_i \arrow[d, "h= \lim_i h_i"] \\
\lim_i T_i \arrow[r]                                                & \lim_i Z_i  \nospacepunct{.}                         
\end{tikzcd}
	\]
	Since $f = h \circ g$ is an effective epimorphism, it follows like in the proof of ii) that $h$ is an equivalence. 
	Thus $h'$ is an equivalence.
	Since ${\mathcal{X}_0}$ is an $\infty$-pretopos, it follows that $g'$ is an inverse limit of effective epimorphisms and so the claim follows from Lemma~\ref{Lemma: limits of eff epis}.
\end{beweis}

\begin{lemma}
	Let ${\mathcal{X}_0}$ be a tiny $\infty$-pretopos.
	Then finite coproducts are universal in $\Pro({\mathcal{X}_0})$.
\end{lemma}
\begin{beweis}
	Again we may use \cite[Proposition 5.3.5.15]{lurie2009higher} to reduce to the case where we are given a cofiltered family of diagrams
	\[
		\begin{tikzcd}
                          & X_i \amalg X_i' \arrow[d, "f_i \ f_i'"] \\
Y_i \arrow[r, "\gamma_i"] & Z_i                                         
\end{tikzcd}
	\]
	and have to show that the induced map
	\[
		 \Big(\lim_i Y_i \times_{\lim_i Z_i} \lim_i X_i\Big) \amalg \Big(\lim_i Y_i \times_{\lim_i Z_i} \lim_i X_i'\Big) \longrightarrow  \lim_i Y_i \times_{\lim_i Z_i} \lim_i\Big(X_i \coprod X_i'\Big)
	\]
	is an equivalence.
	But this map can be identified with the limit of the induced morphisms
	\[
		(Y_i \times_{Z_i} X_i) \amalg (Y_i \times_{Z_i} X_i') \longrightarrow Y_i \times_{Z_i} (X_i \amalg X_i')
	\]
	which are equivalences, as $\mathcal{X}_0$ is an $\infty$-pretopos.
\end{beweis}

Finally we observe that, as a consequence of \cite[Proposition 5.2.8.6]{lurie2012higher} and Proposition~\ref{Proposition: Eff epis sind super in Pro(-)}, the collection of effective epimorphisms in $\Pro({\mathcal{X}_0})$ is closed under composition and it is clearly closed under finite coproducts.
We may thus apply \cite[Proposition~A.3.2.1]{lurie2016spectral} to get the following:

\begin{korollar}\label{Cor: eff epis give topology}
	Let ${\mathcal{X}_0}$ be a bounded $\infty$-pretopos.
	Define a collection of morphisms $\{C_i \rightarrow D\}_{i \in I}$ in $\Pro(\mathcal{X}_0)$ to be covering if and only if there is a finite subset $J \subseteq I$ such that the induced map
	\[
		\coprod_{j \in J} C_j \longrightarrow D
	\]
	is an effective epimorphism in $\Pro({\mathcal{X}_0})$.
	This defines a topology on $\Pro({\mathcal{X}_0})$.
\end{korollar}

\begin{definition}
	\label{def:Pyknotification}
	Let ${\mathcal{X}_0}$ be a tiny bounded $\infty$-pretopos. 
	We call the topology on~$\Pro({\mathcal{X}_0})$ from Corollary~\ref{Cor: eff epis give topology} the \textit{effective epimorphism topology}. 
	For a bounded coherent~$\delta_0$\dash$\infty$\dash topos $\mathcal{X}$, we define the \textit{pyknotification} of $\mathcal{X}$ to be the $\infty$-topos
	\[
		\mathcal{X}^{\pyk} = \Sh_{\eff}^{\hyp}\big(\Pro(\mathcal{X}_{<\infty}^{\operatorname{coh}})\big).
\]	 
\end{definition}

\begin{bemerkung}
	The pyknotification of a bounded coherent $\infty$-topos appeared first in \cite[Construction 3.3.2]{barwick2019pyknotic} under the name \emph{solidifcation}.
	Since the word solidification is also used in \cite{clausen2019condensed} in an unrelated way, a different name is used here.
\end{bemerkung}

\begin{beispiel}
	In the case where $ \mathcal{X} = \mathcal{S} $, the $ \infty $-category of bounded coherent objects in $ \mathcal{S} $ is the $ \infty $-category of $ \pi $-finite spaces $ \mathcal{S}_{\pi} $.
	Let us write $ \mathcal{S}_{\pi}^{^\wedge} $ for the $ \infty $-category $ \Pro(\mathcal{S}_{\pi}) $ of profinite spaces.
	It is shown in \cite[Proposition 13.4.9]{barwick2018exodromy} that any pro-finite set admits and effective epimorphism from a profinite set. 
	In other words the full subcategory $ \Pro(\operatorname{Set}^{\operatorname{fin}}) \subseteq \mathcal{S}_{\pi}^\wedge $ is a basis for the effective epimorphism topology.
	It follows that we have an equivalence of~$\infty$-topoi
	\[
		\mathcal{S}^{\pyk} \simeq \Pyk(\mathcal{S}).
	\]
\end{beispiel}

\begin{beispiel}
	More generally for a qcqs scheme $ X $, any object in $ \Pro(X_{\et}^{\operatorname{constr}}) $ admits an effective epimorphism from an object in $ \Pro(X^{\operatorname{constr}}_{\leq 0}) $ by \cite[Proposition 3.3.8]{barwick2019pyknotic}.
	It follows that $ X_{\et}^{\pyk} $ is the hypercomplete $ \infty $-topos associated to the $ 1 $-topos $ \Sh_{\eff}(\Pro(X^{\operatorname{constr}}_{\leq 0}),\operatorname{Set}) $ introduced in \cite[\S 7.1]{lurie2018ultracategories}. 
	Thus \cite[Example 7.1.7]{lurie2018ultracategories} shows that  $ X_{\et}^{\pyk} $ is equivalent to the hypercomplete $ \infty $-topos $ X_{\proet}^{\hyp} $ of pro-\'etale sheaves on $ X $ defined by Bhatt-Scholze in \cite{bhatt2013pro}.
\end{beispiel}

\begin{bemerkung}
	In principle we can also consider a version of the pyknotification where we consider the $ \infty $-topos $ \Sh_{\eff}(\Pro(\mathcal{X}_{0})) $ of all sheaves with respect to the epimorphism topology instead of just hypersheaves.
	However there are reasons to prefer the hypercomplete version in Definition~\ref{def:Pyknotification}.
	For example in many cases of interest the $ \infty $-topos $ \mathcal{X}^{\pyk} $ will be postnikov complete and even have a set of compact projective generators (see Theorem~\ref{thm:PyknotIsLocWeakCon}), which makes it much easier to work with.
	This will not hold for the non-hypercomplete version in general \cite[Warning 2.2.2]{barwick2019pyknotic}.
\end{bemerkung}

\section{Continuous Straightening-Unstraightening}

Recall that the $ \infty $-category $ \Pyk(\mathcal{S}) $ of pyknotic spaces is equivalent to the $ \infty $-category $ \Fun^\times(\Proj \op,\mathcal{S}) $ of product preserving presheaves on $ \Proj $.
This encourages the following definition \cite[Definition 2.3.1]{barwick2019pyknotic}:

\begin{definition}
	Let $ \mathcal{C} $ be any $ \infty $-category with finite products.
	A \emph{pyknotic object in $ \mathcal{C} $} is a functor $ \Proj \op \to \mathcal{C} $ that preserves finite products.
	We write $ \Pyk(\mathcal{C}) =\Fun^\times(\Proj \op, \mathcal{C}) $ for the $ \infty $-category of pyknotic objects in $ \mathcal{C} $.
	We will simply refer to pyknotic objects in $ \Cat_{\infty} $ as \emph{pyknotic $ \infty $-categories}.
\end{definition}

\begin{definition}
	We call a pyknotic $\infty$-category $\mathcal{C}$ \textit{small}, if $\mathcal{C}(K)$ is a small $ \infty $-category for every $K \in \Proj$.
\end{definition}

\begin{recollection}
	\label{def:FunctorCat}
	Let $\mathcal{C}$ and $\mathcal{D}$ be pyknotic $\infty$-categories. 
	Recall from \cite[Definition 13.3.16]{barwick2018exodromy} that the~$\infty$-category of \textit{continuous functors} from $\mathcal{C}$ to $\mathcal{D}$ is given by
	\[
		\Fun \cts(\mathcal{C},\mathcal{D}) = \int_K \Fun\big(\mathcal{C}(K),\mathcal{D}(K)\big).
	\]
	By \cite[Proposition 2.3]{glasman2016spectrum} the maximal subgroupoid underlying $ \Fun \cts (\mathcal{C},\mathcal{D}) $ is equivalent to the mapping space $ \map_{\Pyk(\Cat_{\infty})}(\mathcal{C},\mathcal{D}) $. 
	In particular the objects in $ \Fun \cts (\mathcal{C},\mathcal{D}) $ are simply morphisms of pyknotic $ \infty $-categories, so natural transformations $ \mathcal{C}(-) \to \mathcal{D}(-) $.
\end{recollection}

\begin{beispiel}
	The $\infty$-category $\Pyk(\mathcal{S})$ of pyknotic spaces may be equipped with a pyknotic structure given by the functor $\Pyk (\mathcal{S})_{/-}$ that takes an extremally disconnected space $K$ to~${\Pyk(\mathcal{S})_{/K}}$ and a morphism $f \colon K \rightarrow K'$ to the functor
	\[
		f^* \colon {\Pyk(\mathcal{S})_{/K}} \longrightarrow {\Pyk(\mathcal{S})_{/K'}}
	\]
	given by pulling back along $f$.
	We denote this pyknotic $ \infty $-category by $ \IPyk(\mathcal{S}) $.
\end{beispiel}

\begin{recollection}[{\cite[Construction 13.3.10]{barwick2018exodromy}}]
	Let $ \mathcal{C} $ be an $ \infty $-category with finite products.
	Evaluating at the point defines a functor
	\[
		\Gamma_* \colon \Pyk(\mathcal{C}) \to \mathcal{C}.
	\]
	For $ X \in \Pyk(\mathcal{C})$, we will refer to $ \Gamma_*(X) $, as the \emph{underlying object of $ X $}.
	If $ \mathcal{C} $ is presentable, the functor $ \Gamma_* $ admits a left adjoint that we will denote
	\[
		(-)^{\operatorname{disc}} \colon \mathcal{C} \to \Pyk(\mathcal{C}).
	\]
	For $ X \in \mathcal{C} $, we refer to $ X^{\operatorname{disc}} $ as the discrete pyknotic object attached to $ X $.
	Furthermore the pyknotic object $ X^{\operatorname{disc}} $ admits an explicit description: Let $ K \in \Proj$ and say that $ K = \lim_{i \in I} K_i $, where every $ K_i $ is a finite set.
	Then we have
	\[
		X^{\operatorname{disc}}(K) \simeq \colim_{i \in I} X^{K_i}
	\]
	where $ X^{K_i} $ denotes the product $ \prod_{K_i} X$.
\end{recollection}

\begin{definition}
	We will say that a continuous functor 
	\[
		f \colon \mathcal{C} \longrightarrow \mathcal{D}
	\]
	in $\Pyk(\Cat_\infty)$ is a \textit{left fibration}, if for every $K \in \operatorname{Proj}$ the functor
	\[
		\mathcal{C}(K) \longrightarrow \mathcal{D}(K)
	\]
	is a left fibration. 
	We will write $\Lfib^{\operatorname{cts}}(\mathcal{C})$ for the full subcategory of $\Pyk(\Cat_\infty)_{/\mathcal{C}}$ spanned by the left fibrations.
\end{definition}

\begin{beispiel}\label{example: discrete left fibrations}
	Let $\mathcal{C}$ be an ordinary $\infty$-category and $p \colon \mathcal{F} \rightarrow \mathcal{C}$ an ordinary left fibration.
	Then this induces a functor $\mathcal{F}^{\operatorname{disc}} \rightarrow \mathcal{C}^{\operatorname{disc}}$ of discrete pyknotic categories. 
	Let $K \in \Proj$ and say that $ K = \lim_{i \in I} K_i $, where $ K_i $ is a finite set for every $ i \in I $.
	Then the induced functor $\mathcal{F}^{\operatorname{disc}}(K) \rightarrow \mathcal{C}^{\operatorname{disc}}(K) $ is given by
	\[
	 	 \colim_i \mathcal{F}^{K_i} \longrightarrow \colim_i \mathcal{C}^{K_i}
	\]
	which, as a filtered colimit of left fibrations, is equivalent to a left fibration.
	Thus the continuous functor $\mathcal{F}^{\operatorname{disc}} \rightarrow \mathcal{C}^{\operatorname{disc}}$ is a left fibration.
\end{beispiel}

The following example will be heavily used in the fourth section:

\begin{beispiel}
	If $ \mathcal{C} $ is an $ \infty $-category, we write $ \Tw(\mathcal{C}) $ for the \emph{twisted arrow $ \infty $-category of $ \mathcal{C} $} \cite[\href{https://kerodon.net/tag/03JG}{Tag 03JG}]{kerodon}.
	For a pyknotic $\infty$-category $\mathcal{C} \in \Pyk(\Cat_\infty)$, we may consider the \textit{continuous twisted arrow $ \infty $-category} $\operatorname{Tw} \cts (\mathcal{C})$. 
	This is the pyknotic $\infty$-category given by the assignment
	\[
		K \longmapsto \operatorname{Tw} (\mathcal{C}(K)).
	\]
	It comes equipped with a canonical functor
	\[
		\operatorname{Tw} \cts (\mathcal{C}) \longrightarrow \mathcal{C} \op \times \mathcal{C}
	\]
	that is a continuous left fibration by \cite[\href{https://kerodon.net/tag/03JQ}{Tag 03JQ}]{kerodon}. 
\end{beispiel}

\begin{recollection}
	\label{rec:Str/Unstr}
	Recall that, for a small $\infty$-category $\mathcal{C}$, there is a natural equivalence
	\[
		\operatorname{Un} \colon \Fun(\mathcal{C},\Cat_{\infty}) \xlongrightarrow{\simeq} \Cocart(\mathcal{C})
	\]
	between $ \Fun(\mathcal{C},\Cat_{\infty}) $ and the $ \infty $-category $ \Cocart(\mathcal{C}) $ of cocartesian fibrations over $ \mathcal{C} $.
	We will denote the inverse of $ \Un $ by $ \operatorname{St} $.
	For a functor $ f \colon \mathcal{C} \to \Cat_{\infty} $ we call $\operatorname{Un}(f)$ the \emph{unstraightening of $ f $}.
	Furthermore the functor $ \Un $ restricts to an equivalence
	\[
		\Fun(\mathcal{C},\mathcal{S}) \xrightarrow{\simeq} \Lfib(\mathcal{C})
	\]
	between $ \Fun(\mathcal{C},\mathcal{S}) $ and the $ \infty $-category of left fibrations over $ \mathcal{C} $.
	By \cite[\href{https://kerodon.net/tag/028T}{Tag 028T}]{kerodon} this equivalence can be explicitly described as follows.
	A functor $ f \colon \mathcal{C} \to \mathcal{S} $ is sent to the left fibration $ \Un(f) \xrightarrow{p} \mathcal{C}$ determined by the pullback square
	\[\begin{tikzcd}
		{\Un(f)} & {\mathcal{S}_{\ast}} \\
		{\mathcal{C}} & {\mathcal{S}\nospacepunct{.}}
		\arrow[from=1-2, to=2-2,,"\pi"]
		\arrow["f"', from=2-1, to=2-2]
		\arrow[from=1-1, to=2-1,"p"]
		\arrow[from=1-1, to=1-2]
	\end{tikzcd}\]
	Here $ \mathcal{S}_\ast $ denotes the $ \infty $-category of pointed spaces and $ \pi $ is the forgetful functor.
	We call $ \pi $ the \emph{universal left fibration}.
\end{recollection}

\begin{definition}
	We define the pyknotic $ \infty $-category of \emph{pointed pyknotic spaces} to be the pyknotic $ \infty $-category
	\[
	 	\IPyk(\mathcal{S}_{\ast}) \colon \Proj \op \to \Cat_{\infty};\; \; \; K \mapsto \Fun^\times((\Proj_{/K}) \op, \mathcal{S}_{\ast}) \simeq \IPyk(\mathcal{S})(K)_{\ast/}.
	\]
	The functor $ \pi \colon \mathcal{S}_{\ast} \to \mathcal{S}$ induces a continuous left fibration $\pi \cts \colon \IPyk(\mathcal{S})_{\ast} \to \IPyk(\mathcal{S})$ that we call the \emph{universal continuous left fibration}.
\end{definition}

We can now prove a pyknotic analogue of the equivalence $ \Lfib(\mathcal{C}) \xrightarrow{\simeq} \Fun(\mathcal{C},\mathcal{S}) $:

\begin{satz}\label{theorem: continuous left fib are cts copresheaves}
	Let $ \mathcal{C} $ be a small pkynotic $ \infty $-category.
	Then there is a canonical equivalence
	\[
		\Un \cts \colon \Fun^{\operatorname{cts}}(\mathcal{C},\mathbf{Pyk}(\mathcal{S})) \xlongrightarrow{\simeq} \Lfib^{\operatorname{cts}}(\mathcal{C})
	\]
	that may be described as follows.
	A continuous functor $ f \colon \mathcal{C} \to \IPyk(\mathcal{S}) $ is sent to the left fibration $ \Un \cts (f) \to \mathcal{C} $ determined by the pullback square
	\[\begin{tikzcd}
		{\Un \cts(f)} & {\IPyk(\mathcal{S}_{\ast})} \\
		{\mathcal{C}} & {\IPyk(\mathcal{S})\nospacepunct{.}}
		\arrow[from=1-2, to=2-2,"\pi \cts"]
		\arrow["f"', from=2-1, to=2-2]
		\arrow[from=1-1, to=2-1]
		\arrow[from=1-1, to=1-2]
	\end{tikzcd}\]
	The inverse sends a continuous left fibration $\mathcal{F} \rightarrow \mathcal{C}$ to the continuous functor whose component at $f \colon K' \rightarrow K \in \Tw(\Proj)$ is given by the assignment
	\begin{align*}
		\mathcal{C}(K) \times (\Proj_{/K'}) \op & \longrightarrow \mathcal{S}\\
		(x_K, \alpha \colon T \rightarrow K') & \longmapsto \operatorname{St}(\mathcal{F}(T))\big((f \circ \alpha)^*(x_K)\big).
	\end{align*}
\end{satz}

For the proof we will need a few more technical details:

\begin{lemma}\label{fiberwise cocart implies cocart}
	Let $p \colon E \rightarrow C$ and $q \colon B \rightarrow C$ be two cocartesian fibrations and let $f \colon E \rightarrow B$ be a functor over $C$ that preserves cocartesian edges.
	Assume that, for every $c \in C$, the induced functor $f_c \colon E_c \rightarrow B_c$ is equivalent to a left fibration.
	Then $f$ is equivalent to a left fibration.
\end{lemma}
\begin{beweis}
	We learned this argument from Alexander Campbell.
	We consider the associated diagram of marked simplicial sets
	\[
		\begin{tikzcd}
E^\natural \arrow[rd, "p"] \arrow[rr, "f"] &          & B^\natural \arrow[ld, "q"] \\
                                           & C^\sharp &                           
\end{tikzcd}
	\]
	(here we use the the notation from \cite[\S 3.1]{lurie2009higher}).
	By the small object argument, we can find a factorization
	\[
		\begin{tikzcd}
                                 & E' \arrow[rd,"f'"] &            \\
E^\natural \arrow[rr,"f"] \arrow[ru,"\iota"] &               & B^\natural
\end{tikzcd}
	\]
	where $\iota$ is marked left anodyne and $f'$ has the right lifting property with respect to marked left anodyne morphisms.
	Thus, by (the dual of) \cite[Proposition 3.1.1.6]{lurie2009higher}, the underlying map of simplicial sets $f' \colon  E' \rightarrow B$ is in particular an isofibration.
	Furthermore it follows that the composite 
	\[
	q \circ f' \colon E' \longrightarrow C^\sharp 
	\]
	has the right lifting property with respect to marked left anodyne morphisms  and thus the underlying map of simplicial sets $E' \rightarrow C$ is a cocartesian fibration.
	So, by \cite[Remark 3.1.3.4]{lurie2009higher} and \cite[Proposition 3.1.3.5]{lurie2009higher}, we get that $\iota$ induces an equivalence on the underlying $\infty$-categories.
	Thus we may replace $f$ by $f'$ to assume that~$f$ is an isofibration.
	As a consequence, the induced map
	\[
		f_c \colon E_c \longrightarrow B_c
	\]
	is an isofibration which is equivalent to a left fibration and hence a left fibration itself.
	Thus we can apply \cite[Proposition 2.4.2.8]{lurie2009higher} and \cite[Proposition 2.4.2.11]{lurie2009higher} to see that~$f$ is a locally cocartesian fibration.
	
	Since the fibers of $f$ are $\infty$-groupoids, it suffices to see that $f$ is a cocartesian fibration.
	We will now prove that every edge $\gamma$ of $E$ is locally $f$-cocartesian. Then the claim will follow from \cite[\href{https://kerodon.net/tag/01V6}{Tag 01V6}]{kerodon}. 
	So let $\gamma \colon x \rightarrow y$ be an edge in $E$. 
	Since $f$ is a locally cocartesian fibration, we may pick a locally cocartesian lift $\delta$ of $f(\gamma)$. 
	This yields a diagram
	\[
	\begin{tikzcd}
\Lambda^2_0 \arrow[rr, "\tau"] \arrow[d] &                                 & E \arrow[d, "f"] \\
\Delta^2 \arrow[r, "\sigma"] \arrow[rru,dotted]            & \Delta^1 \arrow[r, "f(\gamma)"] & B               
\end{tikzcd}
	\]
	where $\sigma$ is the degeneracy given by collapsing the $1 \rightarrow 2$ edge and $\tau$ is given by the diagram
	\[
		\begin{tikzcd}
                                            & {y'} \\
{x} \arrow[ru, "\delta"] \arrow[r, "\gamma"] & {y}
\end{tikzcd}
	\]
	But now, since $\delta$ is locally $f$-cocartesian, there exists a dotted lift in the above diagram.
	This gives rise to a 2-simplex
	\[
	 \begin{tikzcd}
                                            & y' \arrow[rd, "\varepsilon"] &   \\
x \arrow[rr, "\gamma"] \arrow[ru, "\delta"] &                              & y
\end{tikzcd}
	\]
	in $E$, where $\delta$ is locally $f$-cocartesian and $f(\varepsilon) = \id_{f(y)}$.
	Since the functor 
	\[
		f_{p(y)}: E_{p(y)} \longrightarrow B_{p(y)}
	\] 
	is by assumption equivalent to a left fibration and thus conservative, it follows that $\varepsilon$ is an equivalence. Therefore $\gamma$ is also locally $f$-cocartesian.
\end{beweis}

Let us now recall the following result from \cite[Theorem 1.1]{gepner2015lax}:

\begin{satz}\label{Theorem: Oplax colimit = classifying fib}
	Let $\mathcal{C}$ be an $\infty$-category and let $F \colon \mathcal{C} \rightarrow \Cat_\infty$ be a functor. 
	There is a natural equivalence
	\[
		\colim_{\Tw(\mathcal{C})} F(-) \times \mathcal{C}_{-/} \xlongrightarrow{\simeq} \operatorname{Un}(F)
	\]
	of functors $\Fun(\mathcal{C}, \Cat_\infty) \rightarrow \Cat_\infty$. Here $\operatorname{Un}(F)$ denotes the total space of the cocartesian fibration classifying $F$.
\end{satz}

We will need to slightly improve Theorem~\ref{Theorem: Oplax colimit = classifying fib}:

\begin{korollar}
	\label{rem:slightimprovmentofGHN}
	Let $\mathcal{C}$ be an $\infty$-category and let $F \colon \mathcal{C} \rightarrow \Cat_\infty$ be a functor. 
	There is a natural equivalence
	\[
	\colim_{\Tw(\mathcal{C})} F(-) \times \mathcal{C}_{-/} \xlongrightarrow{\simeq} \operatorname{Un}(F)
	\]
	of cocartesian fibrations over $ \mathcal{C} $.
	In particular for any $ f \colon x \to y  $ in $ \mathcal{C} $, the canonical functor
	\[
		c_f \colon F(x) \times \mathcal{C}_{y/} \to \colim_{\Tw(\mathcal{C})} F(-) \times \mathcal{C}_{-/} \xlongrightarrow{\simeq} \operatorname{Un}
	\]
	perserves cocartesian edges.
\end{korollar}
\begin{beweis}
	For any cocartesian fibration $p \colon X \rightarrow \mathcal{C}$ there is a natural equivalence
	\begin{align*}
		\map_{\CoCart(\mathcal{C})}(\colim_{\Tw(\mathcal{C})} F(-) \times \mathcal{C}_{-/}, X) & \simeq \lim_{\Tw(\mathcal{C})} \map_{\CoCart(\mathcal{C})}(F(-) \times \mathcal{C}_{-/},X) \\
		& \simeq \lim_{\Tw(\mathcal{C})} \map_{\Cat_\infty}(F(-),\Fun_{\CoCart(\mathcal{C})}(\mathcal{C}_{-/},X)) \\
		& \simeq \map_{\Fun(\mathcal{C},\Cat_\infty)}(F,\operatorname{St}(X)).
	\end{align*}
	Here the second equivalence follows by adjunction and the last equivalence follows from \cite[Lemma 9.10]{gepner2015lax} and \cite[Proposition 2.3]{glasman2016spectrum}.
	This proves the first part of the claim.
	The second part follows since the explicit formula in Theorem~\ref{Theorem: Oplax colimit = classifying fib} shows that the forgetful functor $\CoCart(\mathcal{C}) \rightarrow {\Cat_\infty}_{/\mathcal{C}}$ preserves colimits.
\end{beweis}

We also recall the following well-known fact:

\begin{lemma}
	\label{rem:claimaboutmaponfibers}
	Let $q \colon \mathcal{E} \rightarrow \mathcal{D}$ and $p \colon \mathcal{D} \rightarrow \mathcal{C}$ be a cocartesian fibrations and $\alpha \colon \Delta^1 \rightarrow \mathcal{D}$ a $p$-cocartesian edge.
	Then the functor classified by the cocartesian fibration
	\[
		\Delta^1 \times_\mathcal{D} \mathcal{E} \rightarrow \Delta^1
	\]
	is given by taking fibers over $\alpha(0)$ in the commutative square
	\[
	\begin{tikzcd}
\mathcal{E}_{p\alpha(0)} \arrow[d] \arrow[rr,"p(\alpha)_!"] && \mathcal{E}_{p\alpha(1)} \arrow[d] \\
\mathcal{D}_{p\alpha(0)} \arrow[rr,"p(\alpha)_!"]           && \mathcal{D}_{p\alpha(1)} \nospacepunct{.}          
\end{tikzcd}
	\]
\end{lemma}
\begin{beweis}
	The above square may be identified with the composite
	\[
	\begin{tikzcd}
{\mathcal{E}_{p\alpha(0)}=\Map_{\mathcal{C}}^\flat(\Delta^{\{0\}},\mathcal{E}^\natural)} \arrow[d] & {\Map^\flat_{\mathcal{C}}((\Delta^1)^\sharp,\mathcal{E}^\natural)} \arrow[l] \arrow[d] \arrow[r] & {\Map^\flat_{\mathcal{C}}(\Delta^{\{1\}},\mathcal{E}^\natural)} \arrow[d] \\
{\mathcal{D}_{p\alpha(0)} = \Map_{\mathcal{C}}^\flat(\Delta^{\{0\}},\mathcal{D}^\natural)}           & {\Map^\flat_{\mathcal{C}}((\Delta^{1})^\sharp,\mathcal{D}^\natural)} \arrow[l] \arrow[r]           & {\Map^\flat_{\mathcal{C}}(\Delta^{\{1\}},\mathcal{D}^\natural)}          
\end{tikzcd}
	\]
	after choosing sections of the left horziontal maps which are induced by precomposition with $\Delta^{\{0\}} \rightarrow (\Delta^1)^\sharp$ and thus trivial fibrations.
	Here we use the notation from \cite[\S 3.1.3]{lurie2009higher}.
	Since $\alpha$ is p-cocartesian we  get an induced map $(\Delta^1)^\sharp \rightarrow \mathcal{D}^\natural$ of marked simplicial sets which induces a commutative cube
	\[
	\begin{tikzcd}[column sep = small]
\mathcal{E}_{\alpha(0)} \arrow[rd] \arrow[dd]                &                                                                     & {\Map^\flat_{\Delta^1}((\Delta^1)^\sharp,(\Delta^1 \times_\mathcal{D} \mathcal{E})^\natural)} \arrow[rd] \arrow[dd] \arrow[ll] &                                                                             \\
                                                             & {\Map_{\mathcal{C}}^\flat(\Delta^0,\mathcal{E}^\natural)} \arrow[dd] &                                                                                                                 & {\Map^\flat_{\mathcal{C}}((\Delta^1)^\sharp,\mathcal{E}^\natural)} \arrow[dd] \arrow[ll] \\
{*} \arrow[rd] &                                                                     & {\Map^\flat_{\Delta^1}((\Delta^1)^\sharp,(\Delta^1)^\sharp)} \arrow[rd] \arrow[ll]                                               &                                                                             \\
                                                             & {\Map_{\mathcal{C}}^\flat(\Delta^0,\mathcal{C}^\natural)}            &                                                                                                                 & {\Map^\flat_{\mathcal{C}}((\Delta^1)^\sharp,\mathcal{D}^\natural)} \arrow[ll]           
\end{tikzcd}
	\]
	where the left and the right face are cartesian.
	We get an analogous cube for the right square in the top rectangle from which the claim follows after choosing sections of the horizontal trivial fibrations above.
\end{beweis}

The following two observations will be needed in order to prove the explicit description of the functor $ \Un \cts $.

\begin{lemma}
	\label{lemma:RightAdjointOfUn}
	The right adjoint of the functor $ \Un \colon \Fun(\mathcal{C},\Cat_{\infty}) \to \Cat_{\infty} $ sends an $ \infty $-category $ \mathcal{X} $ to the functor $ \Fun(\mathcal{C}_{-/},\mathcal{X}) $.
\end{lemma}
\begin{proof}
	For any functor $ F \colon \mathcal{C} \to \Cat_{\infty} $ there is a natural equivalence
	\[
		\map_{\Cat_{\infty}}(\Un(F),\mathcal{X}) \simeq \lim_{\Tw(\mathcal{C})} \map(F(x) \times \mathcal{C}_{y/},\mathcal{X}) \simeq \map_{\Fun(\mathcal{C},\Cat_{\infty})}(F(-),\Fun(\mathcal{C}_{-/},\mathcal{X})).
	\]
	Here the first equivalence follows from Theorem~\ref{Theorem: Oplax colimit = classifying fib} and the second equivalence follows from \cite[Propsotion 2.3]{glasman2016spectrum}.
\end{proof}

Let us denote the right adjoint of the forgetful functor
\[
	F: \CoCart(\mathcal{C}) \to \Cat_{\infty}
\]
by $ R \colon \Cat_{\infty} \to \CoCart(\mathcal{C})$.

\begin{lemma}
	\label{lemma:PullbacksonR}
	Suppose we are given a left fibration $ f \colon \mathcal{X} \to \mathcal{Y} $ of $ \infty $-categories.
	Then commutative square 
	\[\begin{tikzcd}
		{R(\mathcal{X})} & {\mathcal{X}} \\
		{R(\mathcal{Y})} & {\mathcal{Y}}
		\arrow[from=1-1, to=1-2]
		\arrow[from=1-2, to=2-2]
		\arrow[from=2-1, to=2-2]
		\arrow[from=1-1, to=2-1]
	\end{tikzcd}\]
	induced by the adjunction $ F \dashv R $ is a pullback.
\end{lemma}
\begin{proof}
	The composite $ R(\mathcal{Y}) \times_{\mathcal{Y}} \mathcal{X} \to R(\mathcal{Y}) \to \mathcal{C}$ is again a cocartesian fibration. 
	Since $ R(\mathcal{Y}) \times_{\mathcal{Y}} \mathcal{X} \to R(\mathcal{Y}) $ is a left fibration, an edge in $ R(\mathcal{Y}) \times_{\mathcal{Y}} \mathcal{X} $ is cocartesian if and only if its image in $ R(\mathcal{Y}) $ is cocartesian by \cite[Proposition 2.4.1.3]{lurie2009higher}.
	Thus for any other cocartesian fibration $ p \colon \mathcal{E} \to \mathcal{C} $ we have a pullback square
	\[\begin{tikzcd}
		{\map_{\Cocart(\mathcal{C})}(\mathcal{E},R(\mathcal{Y}) \times_{\mathcal{Y}} \mathcal{X})} & {\map_{\Cat_{\infty}}(\mathcal{E},\mathcal{X})} \\
		{\map_{\Cocart(\mathcal{C})}(\mathcal{E},R(\mathcal{Y}) )} & {\map_{\Cat_{\infty}}(\mathcal{E},\mathcal{Y})}
		\arrow[from=2-1, to=2-2]
		\arrow[from=1-2, to=2-2]
		\arrow[from=1-1, to=2-1]
		\arrow[from=1-1, to=1-2]
	\end{tikzcd}\]
	Under the natural equivalence $ \map_{\Cocart(\mathcal{C})}(\mathcal{E},R(-) ) \simeq  \map_{\Cat_{\infty}}(\mathcal{E},-)$ the lower horizontal map identifies with the identity and the right vertical map with the map induced by composition with $ R(f) $.
	Thus we have a natural equivalence
	\[
	\map_{\Cocart(\mathcal{C})}(\mathcal{E},R(\mathcal{Y}) \times_{\mathcal{Y}} \mathcal{X}) \simeq \map_{\Cocart(\mathcal{C})}(\mathcal{E},R(\mathcal{X}) )
	\]
	which proves the claim.
\end{proof}

The strategy is now to first prove a version of Theorem~\ref{theorem: continuous left fib are cts copresheaves} where we replace $ \Pyk(\mathcal{S}) $ by $ \Fun(\Proj \op,\mathcal{S}) $.
In this case the claim will be a rather direct consequence of Theorem~\ref{Theorem: Oplax colimit = classifying fib}.

\begin{definition}
	Let $ \mathcal{C},\mathcal{D} \colon \Proj \op \to \Cat_{\infty}$ be functors.
	\begin{enumerate}
		\item We denote by $\Lfib^{\operatorname{pre-cts}}(\mathcal{C})$ the full subcategory of $ \Fun(\Proj \op , \Cat_\infty)_{/\mathcal{C}} $ spanned by the objectwise left fibrations.
		\item For $ \mathcal{C},\mathcal{D} \colon  \Proj \op \to \Cat_{\infty} $ we define
		\[
		\Fun^{\operatorname{pre-cts}}(\mathcal{C},\mathcal{D}) = \int_K \Fun(\mathcal{C}(K),\mathcal{D}(K)).
		\]
	\end{enumerate}
\end{definition}

As in Recollection~\ref{def:FunctorCat}, we see that the objects of $ \Fun^{\operatorname{pre-cts}}(\mathcal{C},\mathcal{D}) $ are simply natural transformation $ \mathcal{C}(-) \to \mathcal{D}(-) $.

\begin{proposition}
	\label{prop:PreCtsStraighteningUnstraightening}
	Let $ \mathcal{C} \colon \Proj \op \to \Cat_{\infty} $ be a functor.
	Then there is a canonical equivalence
	\[
	 \Un^{\operatorname{pre-cts}} \colon \Fun^{\operatorname{pre-cts}}(\mathcal{C},\Fun((\Proj_{/-}) \op,\mathcal{S})) \xrightarrow{\simeq} \Lfib^{\operatorname{pre-cts}}(\mathcal{C}).
	\]
	that may be described as follows.
	A natural transformation $ \mathcal{C}(-) \to \Fun((\Proj_{/-})\op),\mathcal{S}) $ is sent to the pullback
	\[\begin{tikzcd}
		{\Un^{\operatorname{pre-cts}}(f)} & {\Fun((\Proj_{/-}) \op,\mathcal{S}_{\ast})} \\
		{\mathcal{C}(-)} & {\Fun((\Proj_{/-}) \op,\mathcal{S})\nospacepunct{.}}
		\arrow[from=1-2, to=2-2,"\pi_*"]
		\arrow["f"', from=2-1, to=2-2]
		\arrow[from=1-1, to=2-1]
		\arrow[from=1-1, to=1-2]
	\end{tikzcd}\]
	The inverse sends $\mathcal{F} \rightarrow \mathcal{C} \in \Lfib^{\operatorname{pre-cts}}(\mathcal{C}) $ to the natural transformation, whose component at $f \colon K' \rightarrow K \in \Tw(\Proj)$ is given by the assignment
	\begin{align*}
		\mathcal{C}(K) \times (\Proj_{/K'}) \op & \longrightarrow \mathcal{S}\\
		(x_K, \alpha \colon T \rightarrow K') & \longmapsto \operatorname{St}(\mathcal{F}(T))\big((f \circ \alpha)^*(x_K)\big).
	\end{align*}
\end{proposition}

\begin{beweis}
	Let $q: \tilde{\mathcal{C}} \rightarrow \Proj \op$ denote the cocartesian fibration classifying $\mathcal{C}$.
	By Theorem~\ref{Theorem: Oplax colimit = classifying fib}, we get an equivalence
	\[
		\int_K \Fun\big(\mathcal{C}(K),\Fun((\Proj_{/K}) \op,\mathcal{S})\big) \simeq \Fun\Big(\colim_{\Tw(\mathcal{C})} \mathcal{C}(-) \times \Proj \op_{-/}, \mathcal{S}\Big) \simeq \Lfib(\tilde{\mathcal{C}}).
	\]
	By straightening-unstraightening $ \Lfib^{\operatorname{pre-cts}}(\mathcal{C})$ is equivalent to the subcategory of $ {\Cat_\infty}_{/\tilde{\mathcal{C}}}$ spanned by all functors 
	\[
		f: \mathcal{F} \longrightarrow \tilde{\mathcal{C}}
	\]
	such that the following are satisfied:
	The composite $q \circ f: \mathcal{F} \rightarrow \Proj \op$ is a cocartesian fibration, the morphism~${f \colon \mathcal{F} \rightarrow \tilde{\mathcal{C}}}$ preserves cocartesian edges and, for every $K \in \Proj \op$, the induced map
	\[
		\mathcal{F}_K \longrightarrow \tilde{\mathcal{C}}_K
	\]
	is equivalent to a left fibration. 
	The morphisms in $\Lfib^{\operatorname{pre-cts}}(\mathcal{C})$ are the functors ${E \rightarrow E'}$ over $\tilde{C}$ such that, after composing with $q$, the functor $E \rightarrow E'$ over $\Proj \op$ preserves cocartesian edges.
	By Lemma~\ref{fiberwise cocart implies cocart}, it now follows that every object in  $\Lfib^{\operatorname{pre-cts}}(\mathcal{C})$ is equivalent to a left fibration over $\tilde{\mathcal{C}}$. 
	Furthermore, for any two left fibrations
	\begin{align*}
		g \colon E \longrightarrow \tilde{\mathcal{C}} && \text{and} && h \colon E' \longrightarrow \tilde{\mathcal{C}},
	\end{align*}
	any functor $f \colon E \rightarrow E'$ over $\tilde{\mathcal{C}}$ lies in $\Lfib^{\operatorname{pre-cts}}(\mathcal{C})$. 
	This follows from \cite[Proposition~2.4.1.3]{lurie2009higher} and \cite[Proposition 2.4.2.4]{lurie2009higher}, since an edge $\gamma$ in $E$ is $f \circ g$-cocartesian if and only if $g(\gamma) = h(f(\gamma))$ is $q$-cocartesian and analogously for $E'$.
	So $\Lfib^{\operatorname{pre-cts}}(\mathcal{C})$ and $\Lfib(\tilde{\mathcal{C}})$ are equivalent full subcategories of ${\Cat_\infty}_{/{\tilde{\mathcal{C}}}}$.
	Thus we get an equivalence of~$\infty$\dash categories
	\[
		\Psi \colon \Lfib^{\operatorname{pre-cts}}(\mathcal{C}) \xlongrightarrow{\simeq} \Fun^{\operatorname{pre-cts}}(\mathcal{C},\Fun((\Proj_{/-}) \op,\mathcal{S})).
	\]
	Furthermore it is easy to check that, on the level of objects, $ \Psi $ agrees with assignment of the inverse given in the theorem.
	The above argument and Recolletion~\ref{rec:Str/Unstr} also show that $ \Un^{\operatorname{pre-cts}} $ sends a natural transformation $f \colon \mathcal{C} \to \Fun(\Proj_{/-},\mathcal{S}) $ corresponding to a functor $F \colon \tilde{\mathcal{C}} \to \mathcal{S}$ to the pre-continuous left fibration $ \Un^{\operatorname{pre-cts}}(f) \to \mathcal{C} $ determined by the pullback $ \tilde{\mathcal{C}} \times_{\mathcal{S}} \mathcal{S}_{*} \to \tilde{\mathcal{C}} $ of the universal left fibration along $ F $.
	By Lemma~\ref{lemma:PullbacksonR}, we have a commutative diagram in $ \Cat_{\infty} $
	\[\begin{tikzcd}
		{\tilde{\mathcal{C}} \times_{\mathcal{S}} \mathcal{S}_{*}} & {R(\mathcal{S}_*)} & {\mathcal{S}_*} \\
		{\tilde{\mathcal{C}}} & {R(\mathcal{S})} & {\mathcal{S}}
		\arrow[from=1-1, to=1-2]
		\arrow[from=1-2, to=2-2]
		\arrow[from=2-1, to=2-2]
		\arrow[from=2-2, to=2-3]
		\arrow[from=1-2, to=1-3]
		\arrow[from=1-3, to=2-3]
		\arrow[from=1-1, to=2-1]
	\end{tikzcd}\]
	where all squares are cartesian.
	Therefore Lemma~\ref{lemma:RightAdjointOfUn} shows that under the straightening/unstraightening equivalence the left pullback square translates to a pullback square
	\[\begin{tikzcd}
		{\Un^{\operatorname{pre-cts}}(f)} & {\Fun((\Proj_{/-}) \op,\mathcal{S}_{*/})} \\
		{\mathcal{C}(-)} & {\Fun((\Proj_{/-}) \op,\mathcal{S})}
		\arrow[from=1-2, to=2-2]
		\arrow["f"', from=2-1, to=2-2]
		\arrow[from=1-1, to=2-1]
		\arrow[from=1-1, to=1-2]
	\end{tikzcd}\]
	as claimed.
	\end{beweis}
	
	\begin{beweisv}{Theorem~\ref{theorem: continuous left fib are cts copresheaves}}
	We at first show that the inverse functor in Proposition~\ref{prop:PreCtsStraighteningUnstraightening} restricts to an equivalence
	\[
	\Un \cts \colon \Fun^{\operatorname{cts}}(\mathcal{C},\mathbf{Pyk}(\mathcal{S})) \xlongrightarrow{\simeq} \Lfib^{\operatorname{cts}}(\mathcal{C}).
	\]
	This means that we have to show that a left fibration 
	\[	
		f \colon \tilde{\mathcal{F}} \rightarrow \tilde{\mathcal{C}}
	\]
	classifies a pyknotic $\infty$-category when composed with $\tilde{\mathcal{C}} \rightarrow \Proj\op$ if and only if for every $\alpha \colon K' \rightarrow K$ in $\Tw(\Proj \op)$ the left fibration given by pulling back $f$ along the canonical map
	\[
		c_\alpha \colon \mathcal{C}(K) \times (\Proj_{/K'}) \op \rightarrow \tilde{\mathcal{C}}
	\]
	classifies a functor $F_\alpha \colon \mathcal{C}(K) \times (\Proj_{/K'})\op \rightarrow \mathcal{S}$ which preserves products in the second variable.
	Since $c_\alpha$ preserves cocartesian edges by Corollary~\ref{rem:slightimprovmentofGHN}, it follows from Lemma~\ref{rem:claimaboutmaponfibers} that for any edge
	\[
		\begin{tikzcd}
T' \arrow[rd, "\gamma'",swap] \arrow[rr] &   & T \arrow[ld, "\gamma"] \\
                                    & K' &                       
\end{tikzcd}
	\]
	in $\Proj_{/K'}$ and fixed $x \in \mathcal{C}(K)$ the induced map
	\[
		F_\alpha(x, T) \longrightarrow F_\alpha(x,T')
	\]
	can be identified with the map which is induced by the commutative diagram
	\[
	\begin{tikzcd}
\mathcal{F}( T) \arrow[d] \arrow[r] & \mathcal{F}(T') \arrow[d] \\
\mathcal{C}(T) \arrow[r]           & \mathcal{C}(T')          
\end{tikzcd}
	\]
	after passing to the fibers over $ \gamma^* \alpha^*(x) \in \mathcal{C}(T)$. 
	This shows the claim, since for any $T_0$ and $T_1$ in $(\Proj_{/K'}) \op$ the map $\mathcal{F}(T_0 \amalg T_1) \rightarrow \mathcal{F}(T_0) \times \mathcal{F}(T_1)$ is an equivalence if and only if it is an equivalence after passing to fibers over any $x \in \mathcal{C}(T_0 \amalg T_1)$ in the square
	\[
	\begin{tikzcd}
\mathcal{F}(T_0 \amalg T_1) \arrow[d] \arrow[r] & \mathcal{F}(T_0) \times \mathcal{F}(T_1) \arrow[d] \\
\mathcal{C}(T_0 \amalg T_1) \arrow[r, "\simeq"] & \mathcal{C}(T_0) \times \mathcal{C}(T_1)\nospacepunct{.}          
\end{tikzcd}
	\]
	The remaining claim about the exipicit description of $ \Un \cts $ is now clear since the square
	\[\begin{tikzcd}
		{\IPyk(\mathcal{S}_*)} & {\Fun((\Proj_{/-})\op,\mathcal{S}_\ast)} \\
		{\IPyk(\mathcal{S})} & {\Fun((\Proj_{/-})\op,\mathcal{S})}
		\arrow[from=1-1, to=2-1]
		\arrow[from=1-1, to=1-2]
		\arrow[from=2-1, to=2-2]
		\arrow[from=1-2, to=2-2]
	\end{tikzcd}\]
	is a pullback in $ \Fun(\Proj \op,\Cat_{\infty}) $.
\end{beweisv}

Since over a space any morphism is equivalent to a left fibration, we obtain the following consequence:

\begin{korollar}
Let $F$ be a pyknotic space. Then there is a canonical equivalence
\[
\mathbf{Pyk}(\mathcal{S})_{/F} \xlongrightarrow{\simeq} \Fun^{\operatorname{cts}}\big(F, \mathbf{Pyk}(\mathcal{S})\big).
\]
\end{korollar}

\begin{bemerkung}
	Let $f \colon \mathcal{C} \rightarrow \mathcal{D}$ be a functor of pyknotic $\infty$-categories. 
	It follows from the natuarlity in Theorem~\ref{Theorem: Oplax colimit = classifying fib}, that the equivalences of Theorem~\ref{theorem: continuous left fib are cts copresheaves} fit into a commutative square
	\[\begin{tikzcd}
		{\Fun \cts \big(\mathcal{D}, \Pyk (\mathcal{S})\big)} & {\Fun \cts \big(\mathcal{C}, \mathbf{Pyk}(\mathcal{S})\big)} \\
		{\Lfib \cts (\mathcal{D})} & {\Lfib \cts (\mathcal{C})}
		\arrow["\simeq", from=1-2, to=2-2]
		\arrow["\simeq", from=1-1, to=2-1]
		\arrow[from=2-1, to=2-2]
		\arrow["{f^*}", from=1-1, to=1-2]
	\end{tikzcd}\]
	where the lower hoziontal functor is given by pulling back a continuous left fibration along the continuous functor $ f $.
\end{bemerkung}

\begin{beispiel}
	Let $\mathcal{C}$ be a pyknotic $\infty$-category.
	Then the continuous functor $\mathcal{C}\op \times \mathcal{C} \rightarrow \IPyk(\mathcal{S})$ classified by the continuous twisted diagonal $\Tw \cts (\mathcal{C}) \rightarrow \mathcal{C} \op \times \mathcal{C}$ may be described as follows: 
	For $K \in \Proj$, the component at $K$ is given by
	\begin{align*}
		\mathcal{C}(K)\op \times \mathcal{C}(K) & \longrightarrow {\Pyk(\mathcal{S})_{/K}} \\
		(x,y) & \longmapsto \big((f \colon T \rightarrow K) \mapsto \map_{\mathcal{C}(T)}(f^*x,f^*y)\big).
	\end{align*}
\end{beispiel}

We will now start collecting some first consequences of Theorem~\ref{theorem: continuous left fib are cts copresheaves}.

\begin{bemerkung}\label{remark: notions of discrete agree}
	Let $\mathcal{C}$ be an ordinary $\infty$-category.
	Then the canonical projection map
	\[
		\pr_* \colon \Fun \cts \big(\mathcal{C}^{\operatorname{disc}}, \Pyk(\mathcal{S})\big)
		\longrightarrow \Fun\big(\mathcal{C}, \Pyk(\mathcal{S})\big)
	\]
	is an equivalence of categories. 
	Thus the adjunction between the underlying space and the discrete functor induces an adjunction
	\[
		\adjunction {\operatorname{Disc}_*} {\Fun(\mathcal{C},\mathcal{S})} { \Fun(\mathcal{C}, \Pyk(\mathcal{S}))} {\operatorname{Un}_*}.
	\]
	Furthermore we observe that we have a canonical commutative diagram
	\[
	\begin{tikzcd}
\Lfib(\mathcal{C}) \arrow[r, "\operatorname{Gr}"]                                     & {\Fun ( \mathcal{C}, \mathcal{S})}                                    \\
\Lfib \cts (\mathcal{C}^{\operatorname{disc}}) \arrow[u, "\operatorname{ev}_*"] \arrow[r, "\simeq"] & {\Fun\big(\mathcal{C},\Pyk(\mathcal{S})\big)} \arrow[u, "\operatorname{Un}_*"]
\end{tikzcd}
	\]
	where the horizontal lower arrow is the equivalence from Theorem~\ref{theorem: continuous left fib are cts copresheaves} and $\operatorname{ev}_*$ is the functor given by evaluating at $*$.
	It follows that we may pass to left adjoints vertically to obtain a commutative diagram
	\[
	\begin{tikzcd}
\Lfib(\mathcal{C}) \arrow[r, "\operatorname{Gr}"] \arrow[d, "L"]   & {\Fun ( \mathcal{C}, \mathcal{S})} \arrow[d, "\operatorname{disc}_*"] \\
\Lfib \cts (\mathcal{C}^{\operatorname{disc}}) \arrow[r, "\simeq"] & {\Fun\big(\mathcal{C},\Pyk(\mathcal{S})\big)} \nospacepunct{.}                                
\end{tikzcd}
	\]
	Furthermore, the left adjoint $L$ admits an explicit description given by
	\[
		(\mathcal{F} \rightarrow \mathcal{C}) \longmapsto \big(\mathcal{F}^{\operatorname{disc}} \rightarrow \mathcal{C}^{\operatorname{disc}}\big)
	\]
	and thus agrees with the construction from Example~\ref{example: discrete left fibrations}. 
\end{bemerkung}

\begin{recollection}
	Recall from \cite[Definition 2.1.5]{barwick2018exodromy} that a \emph{stratified space} is an $ \infty $-category $ \Pi $ together with a conservative functor $ \Pi \to P $  to a partially ordered set $ P $.
	We say that $ \Pi \to P $ is \emph{$ \pi $-finite} if $ P $ is finite, $ \Pi $ has only finitely many objects up to equivalence and all mapping spaces of $ \Pi $ are $ \pi $-finite.
	We write $ \Str_{\pi} $ for the full subcategory of $ \Fun(\Delta^1,\Cat_{\infty}) $ spanned by the $ \pi $-finite stratified spaces.
\end{recollection}

\begin{lemma}\label{leftt fibration stratified}
	Let $q \colon \mathcal{F} \rightarrow \Pi$ be a left fibration with $	\pi$-finite fibers and let $\Pi \rightarrow P$ be a~$\pi$\dash finite stratified space. 
	Then $\mathcal{F} \rightarrow P$ is a $\pi$-finite stratified space.
\end{lemma}
\begin{beweis}
	The only thing that is not obvious is that the mapping space $\map_\mathcal{F}(x,y) $ is $\pi$-finite for all $x,y \in \mathcal{F}$.
	Since $ q $ is a left fibration every morphism in $ \mathcal{F} $ is $ q $-cocartesian.
	Thus for every $f \colon x \rightarrow y$ in $\map_\mathcal{F}(x,y)$, the square 
	\[
		\begin{tikzcd}[row sep= small]
		{\map_{\mathcal{F}_{p(y)}}(y,y)} \arrow[d] \arrow[r,"f^*"] & {\map_\mathcal{F}(x,y)} \arrow[d] \\
		* \arrow[r,"q(f)"]                                 & {\map_\Pi\big(f(x),f(y)\big)}            
		\end{tikzcd}
	\]
	is a pullback in $\mathcal{S}$ by \cite[Proposition 2.4.4.2]{lurie2009higher}. 
	By assumption, the bottom right and the top left corner are $\pi$-finite and thus the claim follows.
\end{beweis}

\begin{lemma}
	\label{Remark: Stone constr embeds fully faithful into Pyk/K}
	Consider the conintuous functor $j \colon \mathcal{S}_{\pi}^{\disc} \to \IPyk(\mathcal{S}) $ corresponding under adjunction to the embedding $(-)^{\disc} \colon \mathcal{S}_{\pi} \to \Pyk(\mathcal{S}) $.
	Then $ j(K) \colon \mathcal{S}_{\pi}^{\disc}(K) \to \Pyk(\mathcal{S})_{/K}$ is fully faithful for every $ K \in \Proj $.
\end{lemma}
\begin{proof}
	Write $ K = \lim_{i \in I} K_i $ as an inverse limit of finite sets.
	Recall that
	\[
		\mathcal{S}_{\pi}^{\disc}(K) \simeq \colim_i \mathcal{S}_\pi^{K_i}.
	\]
	For $ i \in I $, the functor $ j(K)(x) $ sends a $ \pi $-finite space $ x  $ over $ K_i $ to the pyknotic space over $ K $ given by the pullback
	\[\begin{tikzcd}
		{j(K)(x)} & x \\
		K & {K_i \nospacepunct{.}}
		\arrow[from=1-1, to=1-2]
		\arrow[from=1-2, to=2-2]
		\arrow[from=1-1, to=2-1]
		\arrow[from=2-1, to=2-2]
	\end{tikzcd}\]
	In particular $ j(K) $ factors through the full subcategory $ {\mathcal{S}_{\pi}^\wedge}_{/K} $ of $ \Pyk(\mathcal{S})_{/K} $ spanned by the profinite spaces over $ K $.
	Now let $ x,y \in \mathcal{S}_{\pi}^{\disc}(K) $.
	Since $ I $ is filtered we may assume that there is some $ i \in I $ such that $ x,y \in \mathcal{S}_{\pi}^{K_i} $.
	Replacing $ I $ by $ I_{/i} $, we may furthemore assume that $ i $ is the final object.
	For a map $ j \to i $ in $ I $, let us denote the pullback $ x \times_{K_i} K_j$ by $ x_j $ and analogously for $ y $.
	We have to see that the map
	\[
		\colim_j \map_{\mathcal{S}_\pi^{K_j}}(x_j,y_j) \to \map_{{\mathcal{S}_\pi^\wedge}_{/K}}\big(j(K)(x),j(K)(y)\big)
	\]
	induced by $ j_K $ is an equivalence.
	Composing with the projection $ j(K)(y) \to y $ induces an equivalence
	\[
		\map_{{\mathcal{S}_\pi^\wedge}_{/K}}\big(j(K)(x),j(K)(y)\big) \xrightarrow{\simeq} \map_{{\mathcal{S}_\pi^\wedge}_{/K_i}}\big(j(K)(x),y).
	\]
	Analogously composing with the projections $ y_j \to y  $ induces an equivalence
	\[
			\colim_j \map_{\mathcal{S}_\pi^{K_j}}(x_j,y_j) \xrightarrow{\simeq} \colim_j \map_{\mathcal{S}_\pi^{K_i}}(x_j,y)
	\]
	We obtain a commutative square
	\[\begin{tikzcd}
		{	\colim_j \map_{\mathcal{S}_\pi^{K_j}}(x_j,y_j)} & {\map_{{\mathcal{S}_\pi^\wedge}_{/K}}\big(j(K)(x),j(K)(y)\big)} \\
		{	\colim_j \map_{\mathcal{S}_\pi^{K_i}}(x_j,y)} & {\map_{{\mathcal{S}_\pi^\wedge}_{/K_i}}\big(j(K)(x),y)\nospacepunct{.}}
		\arrow[from=1-1, to=1-2]
		\arrow["\simeq", from=1-2, to=2-2]
		\arrow["\simeq"', from=1-1, to=2-1]
		\arrow[from=2-1, to=2-2]
	\end{tikzcd}\]
	But it is clear that the lower horizontal map is an equivalence, since $ j(K)(x) \simeq \lim_j x_j $ in $ {\mathcal{S}_{\pi}^{\wedge}}_{/K_i} $ and $ y $ is cocompact in $ {\mathcal{S}_{\pi}^{\wedge}}_{/K_i} $.
\end{proof}

\begin{bemerkung}
	One can give an alternative prove of Lemma~\ref{Remark: Stone constr embeds fully faithful into Pyk/K} as follows.
	Let $ K $ be an extremally diconnected set.
	By \cite[Proposition 4.4.18]{barwick2018exodromy} the $ \infty $-category $ \mathcal{S}_{\pi}^{\disc}(K) $ is equivalent to the full subcategory $ \operatorname{Lcc}(K) $ of the $ \infty $-topos $ \Sh(K) $ of sheaves on $ K $ spanned by the locally constant constructible sheaves in the sense of \cite[Definition 2.5.1]{lurie2016spectral}.
	Then \cite[Corollary 4.4]{haine2022descent} provides a fully faithful embedding
	\[
		\Sh(K) \hookrightarrow \Pyk(\mathcal{S})_{/K}.
	\]
	Therefore we also have a fully faithful embedding $ \mathcal{S}_{\pi}^{\disc}(K) \hookrightarrow \Pyk(\mathcal{S})_{/K}  $ and one can check that this embedding agrees with the functor $ j(K) $.
\end{bemerkung}

\begin{korollar}
	Let $ \mathcal{C}  $ be a pkynkotic $ \infty $-category.
	Then composition with $ j \colon \mathcal{S}_{\pi}^{\disc} \to \IPyk(\mathcal{S}) $ induces a fully faithful functor
	\[
		j_* \colon \Fun \cts (\mathcal{C},\mathcal{S}_{\pi}^{\disc}) \to \Fun \cts (\mathcal{C},\IPyk(\mathcal{S})).
	\]
\end{korollar}

\begin{proof}
	For $ K \in \Proj $, the functor
	\[
		j(K)_{*} \colon \Fun(\mathcal{C}(K),\mathcal{S}_{\pi}^{\disc}(K)) \to \Fun(\mathcal{C}(K),\Pyk(\mathcal{S})_{/K})
	\]
	is fully faithful by Lemma~\ref{Remark: Stone constr embeds fully faithful into Pyk/K}. 
	Therefore $ j_* $ is fully faithful as an end of fully faithful functors.
\end{proof}

\begin{nix}
	\label{nix:IntroducingNotation}
	Let $ S $ be a spectral topological space and $ \mathcal{X} $ a spectral $ S $-stratitified $ \infty $-topos \cite[Definition 9.2.1]{barwick2018exodromy}.
	Consider the profinite straitified shape $ \shape{S}{\mathcal{X}}\in \Pro(\Str_{\pi}) $ of $ \mathcal{X} $  \cite[Definition 10.1.4]{barwick2018exodromy}.
	Say that $ \shape{S}{\mathcal{X}} \simeq \lim_{i} (\mathcal{C}_i \to S_i) $ where $ \mathcal{C}_{i} \to  S_i $ is a $ \pi $-finite stratified space.
	Recall that, by \cite[Lemma 13.6.1]{barwick2018exodromy}, the canonical functor
	\[
	\colim_{i \in I} \Fun\big( \mathcal{C}_i, \mathcal{S}_\pi\big) \longrightarrow \Fun \cts \big(\shape S {\mathcal{X}}, \mathcal{S}_\pi^{\disc}\big)
	\]
	is an equivalence.
	It follows from Remark~\ref{remark: notions of discrete agree} and Theorem~\ref{theorem: continuous left fib are cts copresheaves} that we get the following explicit description of the essential image of the functor $ j_* $:
\end{nix}

\begin{lemma}
	\label{Lemma:UsedtobetheRemark}
	With notations as in \ref{nix:IntroducingNotation}, consider the fully faithful functor
	\[
		j_* \colon \Fun \cts (\shape S {\mathcal{X}},\mathcal{S}_{\pi}^{\disc}) \to \Fun \cts (\shape S {\mathcal{X}},\IPyk(\mathcal{S})).
	\]
	Then a continuous functor $ F \colon \shape S {\mathcal{X}} \to \IPyk(\mathcal{S}) $ is in the essenital image of $ j_* $ if and only if there is a discrete left fibration $ \mathcal{F}_i \longrightarrow \mathcal{C}_i $ with $ \pi $-finite fibers such that the continuous left fibration
	\[
			\mathcal{F} = \lim_{j \in I_{/i}} \mathcal{F}_i \times_{\mathcal{C}_i} \mathcal{C}_j \simeq \mathcal{F}_i \times_{\mathcal{C}_i} \shape S {\mathcal{X}} \to \shape S {\mathcal{X}}
	\]
	classifies $ F$.
	In particular $ \mathcal{F} \to S $ is a profinite stratified space by Lemma~\ref{leftt fibration stratified}.
\end{lemma}

\begin{recollection}
	Recall from \cite[Definition 2.3.7]{barwick2018exodromy} that an $ \infty $-category $ \mathcal{C} $ is called \emph{layered}, if every endomorphism in $ \mathcal{C} $ is an equivalence.
	An $ \infty $-category $ \mathcal{C} $ is called \emph{$ \pi $-finite} if it has finitely many objects up to equivalence and all its mapping spaces are $ \pi $-finite. 
	We write $ \operatorname{Lay}_{\pi} $ for the full subcategory spanned by the $ \pi $-finite and layered $ \infty $-categories.
	Note that for any $ \pi $-finite stratified space $ \Pi \to P $, $ \Pi $ is layered and $ \pi $-finite.
\end{recollection}

\begin{proposition}\label{Prop: Pro(-) embeds fully faithful}
	The embedding $\Fun \cts( \shape S {\mathcal{X}},\mathcal{S}_\pi^{\disc}) \rightarrow \Fun \cts (\shape S {\mathcal{X}},\mathbf{Pyk}(\mathcal{S}))$ extends to a fully faithful tiny limit preserving embedding
	\[
		\iota \colon	\Pro\Big(\Fun \cts\big( \shape S {\mathcal{X}},\mathcal{S}_\pi^{\disc}\big)\Big) \hooklongrightarrow \Fun \cts \big(\shape S {\mathcal{X}},\mathbf{Pyk}(\mathcal{S})\big).
	\]
\end{proposition}
\begin{beweis}
	We have to see that, for a cofiltered diagram $Y_\bullet \colon J \rightarrow \Fun \cts( \shape S {\mathcal{X}},\mathcal{S}_\pi^{\disc})$ and and object $X\in \Fun \cts( \shape S {\mathcal{X}},\mathcal{S}_\pi^{\disc})$, the canonical map
	\begin{equation}
		\colim_j \map_{\Fun \cts\big( \shape S {\mathcal{X}},\mathcal{S}_\pi^{\disc}\big)}(Y_j,X) \longrightarrow \map_{\Fun \cts\big( \shape S {\mathcal{X}},\mathbf{Pyk}(S)\big)}\big(\lim_j Y_j, X\big) \label{cocompact check}
	\end{equation}
	is an equivalence. Using Lemma~\ref{Lemma:UsedtobetheRemark} and Theorem~\ref{theorem: continuous left fib are cts copresheaves}, we see that the embedding 
	\[
		\Fun \cts\big( \shape S {\mathcal{X}},\mathcal{S}_\pi\big) \hooklongrightarrow  \Pyk(\Cat_\infty)_{/ \shape S {\mathcal{X}}} 
	\]
	factors through the full subcategory $\Pro(\operatorname{Lay}_\pi)_{/\shape S {\mathcal{X}}} \subseteq \Pyk(\Cat_\infty)_{/ \shape S {\mathcal{X}}}$ (see \cite[Proposition 13.5.2]{barwick2018exodromy}).
	Again by Lemma~\ref{Lemma:UsedtobetheRemark}, there is a map of stratified spaces
	\[
	\begin{tikzcd}
\shape S {\mathcal{X}} \arrow[d] \arrow[r] & \mathcal{C}_i \arrow[d] \\
S \arrow[r]                                & P_i                    
\end{tikzcd}
	\]
	and a discrete left fibration~$X_i \rightarrow \mathcal{C}_i$ with $\pi$-finite fibers such that the map \eqref{cocompact check} is given by the canonical map
	\[ \arraycolsep = 1.2pt
	\begin{array}{rrl}
		\colim_j \map_{\Fun \cts( \shape S {\mathcal{X}},\mathcal{S}_\pi^{\disc})}(Y_j,X)
			& \multicolumn{1}{l}{\simeq}
				&
		\\[7pt]
		\colim_j \map_{\Pro(\operatorname{Lay}_\pi)_{/\mathcal{C}_i}}(Y_j, X_i)
			& \longrightarrow
				& \map_{\Pro(\operatorname{Lay}_\pi)_{/\mathcal{C}_i}}(\lim_j Y_j,X_i)
		\\[7pt]
			& \simeq
				& \map_{\Fun \cts( \shape S {\mathcal{X}},\mathbf{Pyk}(S))}(\lim_j Y_j, X),
	\end{array}\]
	which is an equivalence as $X_i$ and $\mathcal{C}_i$ are in $\operatorname{Lay}_\pi$.
\end{beweis}

\section{The Proof of Theorem~\ref{Theorem: Cts Exodromy}}

The strategy of the proof of Theorem~\ref{Theorem: Cts Exodromy} is to apply \cite[Proposition A.3.4.2]{lurie2016spectral}:

\begin{satz}\label{nices Sag theorem}
	Let $\mathcal{T}$ be a hypercomplete $\infty$-topos. Let $\mathcal{C} \subseteq \mathcal{T}$ be a full subcategory satisfying the following:
		\begin{enumerate}
			\item The $\infty$-category $\mathcal{C}$ is essentially small.
			\item All objects in $\mathcal{C}$ are coherent.
			\item The subcategory $\mathcal{C} \subseteq \mathcal{T}$ is closed under finite coproducts and under fiber products.
			\item Every object in $\mathcal{T}$ admits a cover by objects in $\mathcal{C}$.
		\end{enumerate}
	Then the composite
	\[
		\mathcal{T} \xlongrightarrow{h} \Fun(\mathcal{T}\op,\mathcal{S}) \xlongrightarrow{\operatorname{restrict}} \Psh(\mathcal{C})
	\]
	induces an equivalence of $\infty$-topoi
	\[
		\mathcal{T} \longrightarrow \Sh_{\operatorname{can}}^{\hyp}(\mathcal{C}).
	\]
	Here $\operatorname{can}$ denotes the topology given by declaring a family $\{U_i \rightarrow X\}_i$ in $\mathcal{C}$ to be covering if there is a finite subset $J \subseteq I$ such that the induced morphism
	\[
		\coprod_{j \in J}  U_j \longrightarrow X
	\]
	is an effective epimorphism in $\mathcal{T}$.
\end{satz}

\begin{nix}
	We will now show that the full subcategory 
	\[
		\iota\Pro\Big(\Fun \cts\big(\shape S {\mathcal{X}},\mathcal{S}_\pi\big)\Big) \subseteq \Fun \cts \big(\shape S {\mathcal{X}},\mathbf{Pyk}(\mathcal{S})\big)
	\]
	satisfies the conditions above.
	Note that i)  and iii) are obvious.
\end{nix}

To see condition ii) and iv), we make the following useful observation:

\begin{lemma}\label{Lemma: Precomp. pre. lim and colim}
	Let $f \colon C \rightarrow D$ be a morphism of small pyknotic $\infty$-categories.
	Then the functor
	\[
		f^* \colon \Fun \cts(\mathcal{D},\IPyk(\mathcal{S})) \rightarrow \Fun \cts (\mathcal{C},\IPyk(\mathcal{S}))
	\]
	has both a left and a right adjoint.
\end{lemma}
\begin{beweis}
	For any morphism $\alpha \colon  K' \rightarrow K$ in $\Proj$ the functor
		\[
			\alpha^* \colon {\Pyk(\mathcal{S})_{/K}} \rightarrow {\Pyk(\mathcal{S})_{/K'}}
		\]
		has both a left and a right adjoint and so does for any $T \in \Proj$ the functor
		\[
			\Fun(\mathcal{C}(K),\Pyk(\mathcal{S})_{/T}) \rightarrow \Fun(\mathcal{C}(K'),\Pyk(\mathcal{S})_{/T})
		\]
		given by precomposing with the functor $\mathcal{C}(\alpha)$.
		Of course the analogous statements hold for $\mathcal{D}$ in place of $\mathcal{C}$.
		It follows that the diagram indexing the limit
		\[
			\int_K \Fun(\mathcal{C}(K),\Pyk(\mathcal{S})_{/K})
		\]
		factors through both inclusions $\mathscr{P}r^L \rightarrow \Cat_\infty$ and $\mathscr{P}r^R \rightarrow \Cat_\infty$ (and similar for $\mathcal{D}$).
		Furthermore the functor
		\[
			f_K^* \colon \Fun(\mathcal{D}(K),\Pyk(\mathcal{S})_{/T}) \rightarrow \Fun(\mathcal{C}(K),\Pyk(\mathcal{S})_{/T})
		\]
		also has a left and a right adjoint for any $K,T \in \Proj$.
		Since the inclusions $\mathscr{P}r^L \rightarrow \Cat_\infty$ and $\mathscr{P}r^R \rightarrow \Cat_\infty$ preserve limits by \cite[Proposition 5.5.3.13]{lurie2009higher} and \cite[Proposition 5.5.3.18]{lurie2009higher}, it follows that the induced map on limits
		\[
			f^* \colon \Fun \cts (\mathcal{D},\IPyk(\mathcal{S})) \rightarrow \Fun \cts (\mathcal{C},\IPyk(\mathcal{S}))
		\]
		is a morphism in both $\mathscr{P}r^L$ and $\mathscr{P}r^R$, as desired.
\end{beweis}

\begin{nix}\label{nix: first appearance of U}
	Let us write $\shape S {\mathcal{X}}^\simeq$ for the maximal pyknotic  sub $\infty$-groupoid of $\shape S {\mathcal{X}}$. Then precomposition with the canonical functor
	\[
	\shape S {\mathcal{X}}^\simeq \longrightarrow \shape S {\mathcal{X}}
	\]
	induces a functor
	\begin{alignat*}{3}
		U \colon \Fun \cts \big(\shape S {\mathcal{X}},\mathbf{Pyk}(\mathcal{S})\big) \longrightarrow & \Fun \cts \big(\shape S {\mathcal{X}}^\simeq,\mathbf{Pyk}(\mathcal{S})\big)\\ \simeq & \Pyk(\mathcal{S})_{/\shape S {\mathcal{X}}^\simeq}
	\end{alignat*}
		which preserves all limits and colimits by Lemma~\ref{Lemma: Precomp. pre. lim and colim}.
		In particular it commutes with the truncation functors $ \tau_{\leq n} $ by \cite[Proposition 5.5.6.28]{lurie2009higher}.
		For every $ K \in \Proj$, the induced functor
		\[
			\Fun(\shape S {\mathcal{X}}(K), \Pyk(\mathcal{S})_{/K}) \to 	\Fun(\shape S {\mathcal{X}}^\simeq (K), \Pyk(\mathcal{S})_{/K})
		\]
		is conservative, it follows that $ U $ is conservative as an end of conservative functors.
		As $\Pyk(\mathcal{S})_{/\shape S {\mathcal{X}}^\simeq}$	is a postnikov complete $\infty$-topos, it is now easy to check using \cite[Proposition 5.5.6.26]{lurie2009higher} that the presentable~$\infty$-category~$\Fun \cts \big(\shape S {\mathcal{X}},\mathbf{Pyk}(\mathcal{S})\big)$ is a postnikov complete $\infty$-topos as well.
		In particular it is hypercomplete.
\end{nix}

\begin{bemerkung} 
	We will now show that iv) of Theorem~\ref{nices Sag theorem} is satisfied.
	The rough idea is the following: Since the $ \infty $-category $\Fun \cts( \shape S {\mathcal{X}}, \Pyk(\mathcal{S}))$ is internally to pyknotic spaces an $ \infty $-category of presheaves, every object should be an \textit{internal colimit} of representables indexed by an \emph{internal diagram}.
	In particular if we restrict the indexing diagram to its underlying pyknotic set, every object $\mathcal{F}$ should receive an effective epimorphism from a coproduct of representables indexed by a pyknotic set.
	Furthermore since the pyknotic $ \infty $-category $\shape S {\mathcal{X}}$ is profinite, the representable presheaves will have values in $\Pro(\Fun \cts (\shape S {\mathcal{X}},\mathcal{S}_\pi))$.
	This means that for any $\mathcal{F} \in \Fun \cts( \shape S {\mathcal{X}}, \Pyk(\mathcal{S}))$ we have an effective epimorphism
	\[
		 \coprod_{k \in K} F_k \twoheadrightarrow \mathcal{F} 
	\]
	where $K$ is some pyknotic set and $F_k$ is in $\Pro(\Fun \cts (\shape S {\mathcal{X}},\mathcal{S}_\pi))$.
	However $\coprod_{k \in K} F_k$ might not necessarily be in $\Pro(\Fun \cts (\shape S {\mathcal{X}},\mathcal{S}_\pi))$.
	To fix this we can futher cover the pyknotic set $K$ by profinite sets $K_\alpha$, so that we get an effective epimorphism
\[
	\coprod_\alpha \bigg( \coprod_{k \in K_\alpha} F_k \bigg) \twoheadrightarrow \mathcal{F}.
\]
	It remains to see that $\coprod_{k \in K_\alpha} F_k$ is contained in $\Pro(\Fun \cts (\shape S {\mathcal{X}},\mathcal{S}_\pi))$ if $K_\alpha$ is a profinite set, which we will see in Proposition~\ref{prop:coproductsarestillfinite}.
	Since we do not have enough pyknotic higher category theory at hand yet to make all of this precise, this does not constitute a full proof. 
	Developing pyknotic higher category theory more generally will be the subject of our future work.
	However for our purposes we will get away with using coproducts indexed by pyknotic sets only.
	This simplifies things as the coproduct functor can simply be understood as the functor
	\[
		\Lfib \cts ( K \times \shape S {\mathcal{X}}) \rightarrow \Lfib \cts ( \shape S {\mathcal{X}})
	\]
	given by composing with the projection $K \times \shape S {\mathcal{X}} \rightarrow \shape S {\mathcal{X}}$.
	We will now execute this strategy:
\end{bemerkung}

\begin{nix}\label{twisted arrow comparison map}
	Let $\mathcal{C}$ be a pyknotic $ \infty $-category and let $\mathcal{F} \rightarrow \mathcal{C}$ be a left fibration.
	We define~$\Tw \cts (\mathcal{F},\mathcal{C})$ to be the pullback in the diagram
	\[
	\begin{tikzcd}
\Tw \cts (\mathcal{F},\mathcal{C}) \arrow[d] \arrow[r] & \Tw \cts (\mathcal{C}) \arrow[d]   \\
\mathcal{F} \op \times \mathcal{C} \arrow[r]          & \mathcal{C} \op \times \mathcal{C}\nospacepunct{.}
\end{tikzcd}
	\]	
	Then the commutative square
	\[
	\begin{tikzcd}
\Tw \cts (\mathcal{F}) \arrow[d] \arrow[r]   & \Tw \cts (\mathcal{C}) \arrow[d]   \\
\mathcal{F} \op \times \mathcal{F} \arrow[r] & \mathcal{C} \op \times \mathcal{C}
\end{tikzcd}
	\]
	induces a morphism $\varphi \colon \Tw \cts (\mathcal{F}) \rightarrow \Tw \cts (\mathcal{F},\mathcal{C})$ in $\Lfib \cts(\mathcal{F} \op \times \mathcal{C})$.
\end{nix}

\begin{lemma}\label{twisted arrow defs die selben}
	With the notation from above, the map $\varphi$ is an equivalence.
\end{lemma}
\begin{beweis}
	It suffices to check this levelwise for every $K \in \operatorname{Proj}$. 
	Thus the claim immediately follow from the next lemma.
\end{beweis}

\begin{lemma}
	Let  $ F \to \mathcal{D} $ be a left fibration of $ \infty $-categories, then the map $ \Tw(F) \to \Tw(F,\mathcal{D}) $ defined as in \ref{twisted arrow comparison map} is an equivalence.
\end{lemma}
\begin{beweis}
	Let $x \in F$ be an object.
	Then, by construction, we get pullback squares
	\[
	\begin{tikzcd}
		\mathcal{D}_{f(x)/} \arrow[r] \arrow[d]   & \Tw(F,\mathcal{D}) \arrow[d]          \\
		\mathcal{D} \arrow[r, "\{x \}\times \id"] & F \op \times \mathcal{D}
	\end{tikzcd}
	\]
	and
	\[
	\begin{tikzcd}
		F_{x/} \arrow[r] \arrow[d]                & \Tw (F) \arrow[d]        \\
		F \arrow[d] \arrow[r, "\{x\} \times \id"] & F \op \times F \arrow[d] \\
		\mathcal{D} \arrow[r, "\{x\} \times \id"] & F \op \times \mathcal{D} \nospacepunct{.}
	\end{tikzcd}
	\]
	So we observe that pulling back the comparison map $\alpha \colon \Tw (F) \rightarrow \Tw(F,\mathcal{D})$ along the inclusion 
	\[
	\mathcal{D} \xlongrightarrow{\{x\} \times \id } F \op \times \mathcal{D}
	\]
	is the functor
	\[
	F_{x/} \longrightarrow \mathcal{D}_{f(x)/}
	\]
	induced by $f$, which is an equivalence because $f$ is a left fibration.
	In particular it follows that $\alpha$ is fiberwise an equivalence and thus an equivalence.
\end{beweis}

\begin{lemma}\label{characterizationf of eff epis}
	Let $f \colon \mathcal{F}\rightarrow \mathcal{G}$ be a morphism in $\Lfib \cts (\mathcal{C)}$. 
	Then $f$ is an effective epimorphism in the $\infty$-topos  $\Lfib \cts (\mathcal{C)}$ if and only if $f(K)$ is essentially surjective for every $K \in \Proj$.
\end{lemma}
\begin{beweis}
	Let us denote by $\mathcal{C}^\simeq$ the maximal pyknotic sub-$\infty$-groupoid. 
	Then, as in \ref{nix: first appearance of U}, the inclusion 
	\[
		i \colon \mathcal{C}^\simeq \hooklongrightarrow \mathcal{C}
	\]
	induces a conservative, limit and colimit preserving functor
	\[
		U \colon \Lfib \cts (\mathcal{C}) \longrightarrow \Pyk 				(\mathcal{S})_{/\mathcal{C}^\simeq}
	\]
	given by pulling back along $i$. This functor thus detects effective epimorphisms. 
	Recall that, for an ordinary left fibration $A \rightarrow B$ of $\infty$-categories, the square
	\[
		\begin{tikzcd}
		A^\simeq \arrow[d] \arrow[r] & A \arrow[d] \\
		B^\simeq \arrow[r]           & B          
		\end{tikzcd}
	\]
	is a pullback square. 
	It follows that the induced map $U(f)$ is given by
	\[
		f^\simeq  \colon \mathcal{F}^\simeq \longrightarrow \mathcal{G}^\simeq
	\]
	over $\mathcal{C}^\simeq$. 
	We note that the morphism $f^\simeq$ is an effective epimorphism in $\Pyk(\mathcal{S})_{/\mathcal{C}^\simeq}$ if and only if it is an effective epimorphism in $\Pyk (\mathcal{S})$ after applying the forgetful functor 
	\[
		\Pyk(\mathcal{S})_{/ \mathcal{C}^\simeq} \longrightarrow \Pyk(\mathcal{S}).
	\]
	So $f$ is an effective epimorphism if and only if $f^\simeq$ is an effective epimorphism considered as a map in $\Pyk(\mathcal{S})$.
	Since the inclusion
	\[
		\Pyk(\mathcal{S}) \hooklongrightarrow \Fun(\Proj \op, \mathcal{S})
	\]
	preserves sifted colimits and finite limits, it follows that $f^\simeq$ is an effective epimorphism if and only if $f^\simeq(K)$ is an effective epimorphism for all $K$. But this happens if and only if $f(K)$ is essentially surjective for all $K$.
\end{beweis}

\begin{nix}
	For the sake of brevity, let us write $\mathcal{C} = \shape S {\mathcal{X}}.$
	Now let $\mathcal{F} \rightarrow \mathcal{C}$ be an object in~$\Lfib \cts (\mathcal{C})$. 
	As in \ref{twisted arrow comparison map}, we consider the commutative diagram
	\[
		\begin{tikzcd}
\Tw \cts (\mathcal{F}) \arrow[rr] \arrow[rd] &                                    & \mathcal{F} \op \times \mathcal{F} \arrow[ld] \\
                                             & \mathcal{F} \op \times \mathcal{C} \nospacepunct{.} &                                              
\end{tikzcd}
	\]
	Let us now pick a collection of profinite sets $\{K_\alpha\}_{\alpha \in A}$ with an effective epimorphism 
	\[
		\psi \colon \coprod_\alpha K_\alpha \longrightarrow \left(\mathcal{F} \op \right)^\simeq.
	\]
	Then pulling back the above triangle along $\coprod_\alpha K_\alpha \times \mathcal{C} \rightarrow \left( \mathcal{F} \op \right)^\simeq \times \mathcal{C} \rightarrow \mathcal{F}\op \times \mathcal{C}$ leads to the commutative triangle
	\[
		\begin{tikzcd}
			P = \Tw \cts (\mathcal{F}) \times_{\mathcal{F} \op \times \mathcal{C}} \left(\coprod_\alpha K_\alpha\right) \times \mathcal{C}
				\arrow[rr]
				\arrow[rd] &                                            & 					\coprod_\alpha K_\alpha \times \mathcal{F} \arrow[ld] \\
			& \coprod_\alpha K_\alpha \times \mathcal{C}\nospacepunct{.} &                                                      
		\end{tikzcd}
	\]
	Let us furthermore write
	\[
		P_\alpha = \Tw \cts (\mathcal{F}) \times_{\mathcal{F} \op \times \mathcal{C}} K_\alpha \times \mathcal{C}
	\]
	such that $P = \coprod_\alpha P_\alpha$. By composing with the projections we in particular get induced morphisms
	\begin{align*}
		P \longrightarrow \mathcal{F} && \text{ and } && P_\alpha \longrightarrow \mathcal{F}
	\end{align*}
	of continuous left fibrations over $\mathcal{C}$.
\end{nix}

\begin{proposition}
	With the notations from above, the map $P \rightarrow \mathcal{F}$ is an effective epimorphism in $\Lfib \cts (\mathcal{C})$.
\end{proposition}

\begin{beweis}
	By construction we have a commutative square
	\[
	\begin{tikzcd}
		P & {\coprod_\alpha K_\alpha \times \mathcal{F}} \\
		{\Tw \cts(\mathcal{F})} & {\mathcal{F} \op \times \mathcal{F}} \nospacepunct{.}
		\arrow[from=1-1, to=1-2]
		\arrow[from=1-2, to=2-2]
		\arrow[from=1-1, to=2-1]
		\arrow[from=2-1, to=2-2]
	\end{tikzcd}
\]
	In particular $ P \to \mathcal{F} $ factors as $ P \to \Tw \cts (\mathcal{F}) \to \mathcal{F}$.
	The morphism $\Tw \cts (\mathcal{F}) \rightarrow \mathcal{F}$ is clearly an effective epimorphism by Lemma~\ref{characterizationf of eff epis}.
	Furthermore, since
		\[
		\coprod_\alpha K_\alpha \longrightarrow \mathcal{F} \op
		\]
	is levelwise essentially surjective, the induced map $P \rightarrow \Tw \cts (\mathcal{F})$ is levelwise essentially surjective as well.
\end{beweis}

\begin{proposition}
	\label{prop:coproductsarestillfinite}
	The left fibration $P_\alpha \rightarrow \mathcal{C}$ is contained in $\Pro(\Fun \cts (\mathcal{C},\mathcal{S}_\pi))$.
\end{proposition}
\begin{beweis}
	To simplify notations, we will say that $\mathcal{C} \rightarrow S \simeq (\mathcal{C}_i \rightarrow S_i)_i$ in $\Pro(\Str_\pi )$ where all $ \mathcal{C}_i \to S_i$ are $ \pi $-finite stratified spaces.
	First we observe that, since the twisted arrow construction is compatible with limits, it follows that the canonical map
	\[
		\Tw \cts (\mathcal{C}) \longrightarrow \lim_i \Tw \cts (\mathcal{C}_i) \times_{\mathcal{C}_i \op \times \mathcal{C}_i} \mathcal{C} \op \times  \mathcal{C}
	\]
	is an equivalence. Using Lemma~\ref{twisted arrow defs die selben}, we see that $P_\alpha$ is given by the pullback square
	\[
		\begin{tikzcd}[column sep=small]
			P_\alpha \arrow[d] \arrow[rr]          &                                              & \Tw \cts(\mathcal{C}) \arrow[d]    \\
			K_\alpha \times  \mathcal{C} \arrow[r] & \mathcal{F} \op \times 	\mathcal{C} \arrow[r] & \mathcal{C} \op \times \mathcal{C} \nospacepunct{.}
		\end{tikzcd}
	\]
	Thus it suffices to see that $P_\alpha^i$, given by the pullback square
	\[
		\begin{tikzcd}
			P_\alpha^i \arrow[d] \arrow[r]        & \Tw \cts(\mathcal{C}_i) \arrow[d]      \\
			K_\alpha \times \mathcal{C} \arrow[r] & \mathcal{C}_i \op \times \mathcal{C}_i
		\end{tikzcd}
	\]
	lies in $\Pro(\Fun \cts (\mathcal{C},\mathcal{S}_\pi))$ when composed with the projection $K_\alpha \times  \mathcal{C} \rightarrow \mathcal{C}$.
	For this we observe that, since the twisted arrow functor $\Tw(-)$ commutes with filtered colimits the canonical map $ \Tw(\mathcal{C}_i)^{\disc} \to \Tw \cts (\mathcal{C}_i^{\disc}) $ is an equivalence.
	Thus the functor $\Tw \cts (\mathcal{C}_i) \rightarrow \mathcal{C}_i \op \times \mathcal{C}_i$ is a discrete left fibration with $\pi$-finite fibers.
	Let us say that~$K_\alpha = \{K_\alpha^j\}_{j \in J}$ as a profinite set. 
	Then the map $K_\alpha \rightarrow \mathcal{F} \op \rightarrow \mathcal{C}_i \op$ factors through some~$K_\alpha \rightarrow K_\alpha^{j_0}$. 
	It follows that the induced map 
	\[
		P_\alpha^i \longrightarrow \lim_{j \in J_{/j_0}} \big(K_\alpha^j \times \mathcal{C} \times_{\mathcal{C}_i \op \times \mathcal{C}_i} \Tw \cts(\mathcal{C}_i)\big)
	\]
	is an equivalence.
	Thus it suffices to see that the composite
	\[
		K_\alpha^j \times \mathcal{C} \times_{\mathcal{C}_i \op \times \mathcal{C}_i} \Tw \cts(\mathcal{C}_i) \longrightarrow K_\alpha^j \times  \mathcal{C} \longrightarrow \mathcal{C}
	\]
	is contained $\Pro(\Fun \cts (\mathcal{C},\mathcal{S}_\pi))$ for all $j$. 
	But by construction all squares in the diagram 
	\[
	\begin{tikzcd}[column sep=small]
K_\alpha^j \times \mathcal{C} \times_{\mathcal{C}_i \op \times \mathcal{C}_i} \Tw \cts(\mathcal{C}_i) \arrow[d] \arrow[r] & K_\alpha^j \times \mathcal{C}_i \times_{\mathcal{C}_i \op \times \mathcal{C}_i} \Tw \cts(\mathcal{C}_i) \arrow[d] \\
K_\alpha^j \times \mathcal{C} \arrow[d] \arrow[r]                                                                         & K_\alpha^j \times \mathcal{C}_i \arrow[d]                                                                         \\
\mathcal{C} \arrow[r]                                                                                                     & \mathcal{C}_i                                                                                                    
	\end{tikzcd}
	\]
	are pullback squares and thus the claim follows, as the map
	\[
		K_\alpha^j \times \mathcal{C}_i \times_{\mathcal{C}_i \op \times \mathcal{C}_i} \Tw \cts (\mathcal{C}_i) \longrightarrow K_\alpha^j \times \mathcal{C}_i \longrightarrow \mathcal{C}_i
	\]
	is a discrete left fibration with $\pi$-finite fibers.
\end{beweis}

Combining the last two propositions, we obtain the following:

\begin{korollar}
Every object in the $\infty$-topos $\Fun \cts \big(\shape S {\mathcal{X}}, \IPyk (\mathcal{S})\big)$ admits a cover by objects in $\Pro\big(\Fun \cts \big(\shape S {\mathcal{X}},\mathcal{S}_\pi^{\disc}\big)\big)$.
\end{korollar}

We will now show that ii) of Theorem~\ref{nices Sag theorem} is satisfied:

\begin{proposition}
	The fully faithful embedding
	\[
	 \iota \colon \Pro\Big(\Fun \cts\big( \shape S {\mathcal{X}},\mathcal{S}_\pi^{\disc}\big)\Big) \hooklongrightarrow \Fun \cts \big(\shape S {\mathcal{X}},\mathbf{Pyk}(\mathcal{S})\big)
	\]
	 factors through the full subcategory spanned by the coherent objects.
\end{proposition}
\begin{beweis}
	We will show that all objects in $\Pro(\Fun \cts (\shape S {\mathcal{X}},\mathcal{S}_\pi^{\disc}))$ are $n$-coherent using induction on $n$. Let us start with $n=0$.
	Recall that the functor
	\[
	 	U \colon \Fun \cts \big(\shape S {\mathcal{X}},\mathbf{Pyk}(\mathcal{S})\big) \longrightarrow \Pyk(\mathcal{S})_{/\shape S {\mathcal{X}}^\simeq}
	\]
	is conservative and preserves all limits and colimits by Lemma~\ref{Lemma: Precomp. pre. lim and colim}. 
	Thus it suffices to see that, for an object $\mathcal{F} \in \Pro(\Fun \cts (\shape S {\mathcal{X}},\mathcal{S}_\pi^{\disc}))$, the pyknotic space $U(\mathcal{F})$ over $\shape S {\mathcal{X}}^\simeq$ is quasi-compact. 
	We now observe that the functor $U$ takes objects in $\Fun \cts (\shape S {\mathcal{X}},\mathcal{S}_\pi^{\disc})$ to objects in $\Fun \cts(\shape S {\mathcal{X}}^\simeq ,\mathcal{S}_\pi^{\disc}) \subseteq {\mathcal{S}_\pi^\wedge}_{/\shape S {\mathcal{X}}^\simeq}$.
	Since the inclusion
	\[
		{\mathcal{S}_\pi^\wedge}_{/\shape S {\mathcal{X}}^\simeq} \longrightarrow \Pyk(\mathcal{S})_{/\shape S {\mathcal{X}}^\simeq}
	\]
	preserves limits, it follows that $U$ takes objects in~$\Pro(\Fun \cts (\shape S {\mathcal{X}},\mathcal{S}_\pi^{\disc}))$ to objects in the full subcategory ${\mathcal{S}_\pi^\wedge}_{/\shape S {\mathcal{X}}^\simeq }$. 
	Now \cite[Corollary 13.4.10]{barwick2018exodromy} and \cite[Propositions 3.1.3, Proposition 2.2.2 and Remark 2.0.5]{lurie2016spectral} imply that all object in ${\mathcal{S}_\pi^\wedge}_{/\shape S {\mathcal{X}}^\simeq }$ are coherent and thus in particular quasi-compact. 
	This completes the case $n=0$. 
	The induction step is now clear from \cite[Corollary A.2.1.4]{lurie2016spectral}.
\end{beweis}

To complete the proof of Theorem~\ref{Theorem: Cts Exodromy}, we will need a few more technical details:

\begin{lemma}\label{lemma: effective epis come from finite stages}
	Let $X_\bullet \colon I \rightarrow \mathbf{bPretop}_\infty^{\delta_0}$ be a tiny filtered diagram in the $\infty$-category of tiny bounded $\infty$-pretopoi.
	Let $ X  = \colim_i X_i$ denote the colimit in $ \mathbf{bPretop}_\infty^{\delta_0} $ and let $f \colon C \rightarrow D$ be an effective epimorphism in $ X $.
	Then there is an $i_0 \in I$ and an effective epimorphism $f_{i_0} \colon C_{i_0} \rightarrow D_{i_0}$ mapping to $f$ under the canonical functor
	\[
		k_{i_0}\colon X_{i_0} \longrightarrow X.
	\]
\end{lemma}
\begin{beweis}
	By \cite[Proposition A.8.3.1]{lurie2016spectral}, the inclusion $\mathbf{bPretop}_\infty^{\delta_0} \rightarrow \Cat_\infty^{\delta_0}$ preserves filtered colimits.
	Thus we may find a morphism $f_{j_0} \colon C_{j_0} \rightarrow D_{j_0}$ such that $k_{j_0}(f_{j_0})$ is equivalent to $f$.
	Since $k_{j_0}$ is a morphism of $\infty$-pretopoi, we get an equivalence
	\[
		k_{j_0}(\check{C}(f_{j_0})_\bullet) \simeq \check{C}(f)_\bullet
	\]
	of simplicial objects in $X$.
	Furthermore, since $k_{j_0}$ is a morphism of $\infty$-pretopoi, it preserves finite limits and effective epimorphisms and thus geometric realizations of groupoid objects.
	It follows that the canonical map
	\[
		c \colon \lvert \check{C}(f_{j_0})_\bullet \rvert \longrightarrow D_{j_0}
	\]
	becomes an equivalence after applying $k_{j_0}$.
	Thus there is a map $\gamma \colon j_0 \rightarrow i_0$ such that $X_\gamma(c)$ is an equivalence and since $X_\gamma$ is a morphism of $\infty$-pretopoi, it follows that $X_\gamma(f_{j_0})$ is an effective epimorphism, as desired.
\end{beweis}

\begin{korollar}\label{Cor: inlusion from constr. sheaves preserves eff. epis}
	Let $K= \{K_i\}_i$ be a profinite space. 
	Then the fully faithful functor 
	\[
		\Fun \cts (K, \mathcal{S}_\pi^{\disc}) \longrightarrow {\Pyk(\mathcal{S})_{/K}}
	\]
	preserves effective epimorphisms.
\end{korollar}
\begin{beweis}
	By Lemma~\ref{lemma: effective epis come from finite stages}, it suffices to see that, for all $i$, every effective epimorphism $f$ in~$\Fun(K_i, \mathcal{S}_\pi)$ maps to an effective epimorphism in ${\Pyk(\mathcal{S})_{/K}}$.
	Denote by $ p_i \colon K \to K_i $ the projection.
	Since we have a commutative diagram
	\[
		\begin{tikzcd}
{\Fun(K_i,\mathcal{S}_\pi^{\disc})} \arrow[d,"p_i^*"] \arrow[r,"\varphi"] & {\Pyk(\mathcal{S})_{/K_i}} \arrow[d,"p_i^*"] \\
{\Fun \cts (K,\mathcal{S}_\pi^{\disc})} \arrow[r]       & {\Pyk(\mathcal{S})_{/K}}            
\end{tikzcd}
	\]
	it suffices to see that the top horizontal functor $\varphi$ preserves effective epimorphism.
	Picking a section of the canonical morphism $K_i \rightarrow \pi_0(K_i)$ and precomposing with it, we may assume that $K_i$ is a finite set.
	In this case $\varphi$ is identified with the product of finitely many copies of the inclusion
	\[
		\mathcal{S}_\pi \longrightarrow \Pyk (\mathcal{S}),
	\]
	which clearly preserves effective epimorphisms. This completes the proof.
\end{beweis}

We finally arrive at the following:

\begin{proposition}\label{Proposition: Notions of effective epis agree}
	A  morphism $f \colon X \rightarrow Y$ in the $\infty$-category $\Pro (\Fun \cts(\shape S {\mathcal{X}},\mathcal{S}_\pi)^{\disc})$ is an effective epimorphism if and only if $\iota(f)$ is an effective epimorphism in the $ \infty $-topos $\Fun \cts \big(\shape S {\mathcal{X}},\mathbf{Pyk}(\mathcal{S})\big)$.
\end{proposition}
\begin{beweis}
	Let us first assume that $\iota(f)$ is an effective epimorphism. 
	We then may factor $f \simeq g \circ h$, where $h$ is an effective epimorphism and $g$ is $(-1)$-truncated.
	Since the inclusion $\iota$ preserves finite limits, it follows that $\iota(g)$ is $(-1)$-truncated as well.
	But by \cite[Corollary 6.2.3.12]{lurie2009higher}, the map $\iota(g)$ is an effective epimorphism because $\iota(f)$ is.
	This implies that $\iota(g)$, and thus $g$, is an equivalence, as desired.
	
	Now we show that $\iota$ preserves effective epimorphisms.
	Again we consider the inclusion~$K = \shape S {\mathcal{X}}^\simeq \hookrightarrow \shape S  {\mathcal{X}}$. 
	Pre-composing with this inclusion induces a morphism of $\infty$-pretopoi
	\[
		\Fun \cts\big(\shape S {\mathcal{X}},\mathcal{S}_\pi^{\disc}\big) \longrightarrow \Fun \cts(K,\mathcal{S}_\pi^{\disc}).
	\]
	So it follows from Proposition~\ref{Proposition: Eff epis sind super in Pro(-)} that the induced functor
	\[
		\Pro\Big(\Fun \cts\big(\shape S {\mathcal{X}},\mathcal{S}_\pi^{\disc}\big)\Big) \longrightarrow \Pro\big( \Fun \cts(K,\mathcal{S}_\pi^{\disc})\big)
	\]
	preserves effective epimorphisms.
	We may thus reduce to showing that the induced functor 
	\[
		j \colon \Pro\big(\Fun \cts (K, \mathcal{S}_\pi^{\disc})\big) \longrightarrow {\Pyk(\mathcal{S})_{/K}}
	\]
	preserves effective epimorphisms.
	By Corollary~\ref{Cor: inlusion from constr. sheaves preserves eff. epis}, the inclusion
	\[
		k \colon \Fun \cts (K, \mathcal{S}_\pi^{\disc}) \hooklongrightarrow \Pyk (\mathcal{S})_{/K}.
	\]
		preserves effective epimorphisms. Furthermore $k$ factors through the full subcategory~${\mathcal{S}_\pi^\wedge}_{/K}$ and hence so does $j$.
		By Lemma~\ref{Lemma: limits of eff epis} and since effective epimorphisms in slice categories are detected by the projections, it suffices to see that the inclusion $\mathcal{S}_\pi^\wedge \rightarrow \Pyk(\mathcal{S})$ preserves effective epimorphisms, which is clear by \cite[Corollary 13.4.10]{barwick2018exodromy}.
\end{beweis}

We have finally collected all the necessary ingredients that are needed to prove our main theorem:

\begin{beweisv}{\ref{Theorem: Cts Exodromy}}
	The Exodromy Theorem provides an equivalence of tiny $\infty$-pretopoi
	\[
		\Fun \cts\big(\shape S {\mathcal{X}},\mathcal{S}_\pi^{\disc}\big) \simeq \mathcal{X}_{< \infty}^{\coh}.
	\] 
	We have seen above that the full subcategory 
	\[
		\Pro \Big(\Fun \cts\big(\shape S {\mathcal{X}},\mathcal{S}_\pi^{\disc}\big)\Big) \hooklongrightarrow \Fun^{\operatorname{cts}}\big(\shape S {\mathcal{X}}, \mathbf{Pyk}(\mathcal{S})\big)
	\]
	satisfies the assumptions of Theorem~\ref{nices Sag theorem} and thus we get an equivalence
	\[
		\Sh_{\operatorname{can}}^{\hyp}\Big(\Pro (\Fun \cts\big(\shape S {\mathcal{X}},\mathcal{S}_\pi^{\disc})\big)\Big) \simeq \Fun^{\operatorname{cts}}\big(\shape S {\mathcal{X}}, \mathbf{Pyk}(\mathcal{S})\big).
	\]
	Thus it remains to see that the topologies $\operatorname{can}$ and $\eff$ on $\Pro (\Fun \cts(\shape S {\mathcal{X}},\mathcal{S}_\pi))$ agree, but this is just a reformulation of Proposition~\ref{Proposition: Notions of effective epis agree}.
\end{beweisv}

As a first easy consequence, we obtain the following :

\begin{korollar}
	\label{cor: pullback pre limits}
	Let $f: X \rightarrow Y$ be any morphism of schemes.
	Then the induced pull-back functor
	\[
		f^*: Y_{\proet}^{\hyp} \longrightarrow X_{\proet}^{\hyp}
	\]
	has both a left and a right adjoint.
\end{korollar}
\begin{beweis}
	By the adjoint functor theorem \cite[Corollary 5.5.2.9]{lurie2009higher} it suffices to see that $ f^* $ preserves all limits and colimits.
	We may cover $ X $ by affine opens $j_i \colon U_i \to X $ such that for every $ i $ we have an affine open $t_i \colon V_i \to Y $ and a commutative diagram
	\[
	\begin{tikzcd}
		X & Y \\
		{U_i} & {V_i}\nospacepunct{.}
		\arrow["f", from=1-1, to=1-2]
		\arrow["{f_i}"', from=2-1, to=2-2]
		\arrow["{j_i}",hook, from=2-2, to=1-2]
		\arrow["{t_i}",hook, from=2-1, to=1-1]
	\end{tikzcd}
	\]
	Since the $ j_i ^*$ are jointly conservative and restrictions to open subschemes commute with all limits and colimits, it suffices to see that each $ t_i^* \circ f^* $ preserves limts.
	Because $ t_i^* \circ f^* \simeq f_i^* \circ j_i^* $, it suffices to see that $ f_i^* $ preserves all limits and colimits and we may therefore assume that $ X $ and $ Y $ are affine, so in particular quasi-compact.
	In this case $f^*$ corresponds to the functor
	\[
		\Gal(f)^* \colon \Fun \cts \big( \Gal(Y), \IPyk(\mathcal{S})\big) \longrightarrow  \Fun \cts \big( \Gal(X), \IPyk(\mathcal{S})\big)
	\]
	given by precomposing with $\Gal(f) \colon \Gal(X) \rightarrow \Gal(Y)$ via Corollary~\ref{cor:X proet= Psh(Gal(X)op)}.
	Thus the claim follows from Lemma~\ref{Lemma: Precomp. pre. lim and colim}.
\end{beweis}

\begin{bemerkung}
	We will now roughly sketch how to circumvent the enlargement of universes which appears in our results, following \cite[\S 1.4]{barwick2019pyknotic}.
	Let $\mathcal{X}$ be a spectral $\infty$-topos and let $\beta$ be an uncountable regular cardinal such that $\shape S {\mathcal{X}}$ is a $\beta$-small inverse limit of $\pi$-finite layered $\infty$-categories.
	Let $\Pro(\mathcal{X}^{\coh}_{<\infty})_\beta$ denote the small subcategory spanned by the $\beta$-cocompact objects. 
	We define
	\[
		\mathcal{X}^{\pyk,\beta} = \Sh_{\eff}^{\hyp}(\Pro(\mathcal{X}^{\coh}_{<\infty})_\beta)
	\]
	and observe that $\shape S {\mathcal{X}}$ naturally defines a sheaf of $ \infty $-categories on $\Pro(\mathcal{S}_\pi)_\beta$, i.e. an $  \infty $-category object in $\mathcal{S}^{\pyk,\beta}$.
	Furthermore let us write $ \Pyk(\mathcal{S})^\beta = \mathcal{S}^{\pyk,\beta}$
	We can then reproduce the results of \S 3 and \S 4 in this framework to obtain an equivalence
	\[
		\mathcal{X}^{\pyk,\beta} \simeq \Fun \cts_\beta(\shape S {\mathcal{X}},\IPyk(\mathcal{S})^\beta)
	\]
	of $\infty$-topoi where we write
	\[
		\Fun \cts_\beta(\shape S {\mathcal{X}},\IPyk(\mathcal{S})^\beta) \simeq \int_{T \in (\mathcal{S}_{\pi}^{\wedge})_\beta} \Fun (\shape S {\mathcal{X}}(T), \Pyk(\mathcal{S})^\beta_{/T}).
	\]
	Considering the left Kan-extensions along $\Pro(\mathcal{X}_{<\infty}^{\coh})_{\lambda_0} \hookrightarrow \Pro(\mathcal{X}_{<\infty}^{\coh})_{\lambda_1}$ for any $ \beta < \lambda_0 < \lambda_1$, we obtain an equivalence
	\[
		\mathcal{X}^{\pyk,\operatorname{acc}} \coloneqq  \colim_{\lambda > \beta} \mathcal{X}^{\pyk,\lambda} = \Sh_{\eff}^{\hyp,\operatorname{acc}}(\Pro(\mathcal{X}^{\coh}_{<\infty})) \xrightarrow{\simeq} \colim_{\lambda > \beta} \Fun \cts_\lambda(\shape S {\mathcal{X}},\IPyk(\mathcal{S})^\lambda).
	\]
	Furthermore the filtered colimit
	\[
		\colim_{\lambda > \beta} \Pyk(\mathcal{S})^\lambda = \Sh_{\eff}^{\hyp,\operatorname{acc}}(\mathcal{S}_\pi^\wedge)
	\]
	is given by the $\infty$-category of accessible sheaves on $\mathcal{S}_\pi^\wedge$, which we can further identify with the $\infty$-category of \emph{condensed spaces} $\Cond(\mathcal{S})$ of Clausen and Scholze.
	We may consider $ \Cond(\mathcal{S}) $ as a hypersheaf with respect to the effective epimorphism topology on $ \mathcal{S}_{\pi}^{\wedge} $ as follows.
	Denote by $ \mathbf{Cond}(\mathcal{S})^{\lambda} $ the sheaf given by left Kan-extension of
	\[
		\IPyk(\mathcal{S})^\lambda \colon \Pro(\mathcal{S}_\pi)_\lambda \to \Cat_\infty
	\]
	along $ (\mathcal{S}_{\pi}^{\wedge})_\lambda \to \mathcal{S}_{\pi}^{\wedge} $ and define
	\[
		\mathbf{Cond}(\mathcal{S}) \simeq \colim_{\lambda} \mathbf{Cond}(\mathcal{S})^\lambda.
	\]
	We may also consider $\shape S {\mathcal{X}}$ as a sheaf of $ \infty $-categories on $ \mathcal{S}_{\pi}^{\wedge} $ via left Kan-extension.
	The resulting sheaf is therefore an accessible sheaf and thus $ \kappa $-compact for some regular cardinal $ \kappa $.
	It follows that the $ \infty $-category of natural transformations
	\[
	  \Fun^{\Cond}(\shape S {\mathcal{X}}, \mathbf{Cond}(\mathcal{S})) = \int_{K \in \mathcal{S}_{\pi}^{\wedge}} \Fun(\shape S {\mathcal{X}}(K),\mathbf{Cond}(\mathcal{S})(K))
	\]
	is equivalent to the filtered colimit
	\[
	\colim_{\lambda > \kappa} \Fun \cts_\lambda(\shape S {\mathcal{X}},\IPyk(\mathcal{S})^\lambda).
	\]
	Thus the above equivalence shows that we have an equivalence
	\[
		\mathcal{X}^{\pyk,\operatorname{acc}} \simeq \Fun^{\Cond}(\shape S {\mathcal{X}}, \mathbf{Cond}(\mathcal{S})).
	\]
\end{bemerkung}

\section{Local Contractibility}

We consider the following $\infty$-categorical analogue of \cite[Definition 3.2.1]{bhatt2013pro}:

\begin{definition}
	Let $\mathcal{X}$ be an $\infty$-topos. 
	Then an object $P \in \mathcal{X}$ is called \textit{weakly contractible}, if the functor
	\[
		\map_\mathcal{X}(P,-) \colon \mathcal{X} \longrightarrow \mathcal{S}
	\]
	preserves geometric realizations of simplicial objects.
	We say that $\mathcal{X}$ is \textit{locally weakly contractible}, if every object $X \in \mathcal{X}$ admits a cover by coherent weakly contractible objects.
\end{definition}

\begin{bemerkung}
	\label{rem:WCONTRpreEpis}
	If $ P \in \mathcal{X} $ is locally weakly contractible it follows that $ \map_{\mathcal{X}}(P,-) $ preserves effective epimorphisms.
	As an easy consequence we obtain the following Lemma.
\end{bemerkung}

\begin{lemma}
	\label{Lem:WContrSections}
	Let $ P \in \mathcal{X} $ be weakly contractible. 
	Then every effective epimorphism $p \colon Y \to P $ splits.
\end{lemma}
\begin{beweis}
	By Remark~\ref{rem:WCONTRpreEpis} the induced morphism
	\[
		\map_{\mathcal{X}}(P,Y) \to \map_{\mathcal{X}}(Y,Y)
	\]
	is an effective epimorphism in $ \mathcal{S} $, so surjective on connected components.
	In particular $ \id_{Y} $ admits a preimage $ s \colon Y \to P $ under $ p_* $, as desired.
\end{beweis}

\begin{lemma}\label{Lemma: retracts of coh are coh}
	Let $\mathcal{X}$ be a locally $ n $-coherent $\infty$-topos.
	Then $ (n+1) $-coherent objects in $\mathcal{X}$ are stable under retracts.
	In particular coherent objects are stable under retracts in a locally coherent $ \infty $-topos.
\end{lemma}
\begin{beweis}
	Let $ k \geq n +1 $.
	We will show that retracts of $ k $-coherent objects are $k$-coherent using induction over~$k$.
	The case $k=0$ is clear, since retracts of quasi-compact objects are quasi-compact.
	So let $ 0 < k \leq n $ and assume that retracts of $ k $-coherent objects are $k$-coherent.
	Let $X$ be a retract of a $ (k+1) $-coherent object $X'$. 
	Since $\mathcal{X}$ is locally $ k $-coherent, it suffices by \cite[Corollary A.2.1.4]{lurie2016spectral} to see that for every cospan
	\[
		U \longrightarrow X \longleftarrow V
	\]
	where $U$ and $V$ are $ k $-coherent, the pullback $U \times_X V$ is $k$-coherent.
	But since $X$ is a retract of~$X'$, the pullback $U \times_X V$ is a retract of $U \times_{X'} V$, which is $ k $-coherent as $ X' $ is $ (k+1) $-coherent.
	So $U \times_X V$ is $k$-coherent by assumption and thus $X$ is $(k+1)$-coherent, as desired.
\end{beweis}

\begin{construction}
	Let $ \mathcal{X} $ be a spectral $ \infty $-topos. Let $ K = \shape S {\mathcal{X}}^\simeq$.
	We again consider the induced functor
		\[
		U \colon \Fun \cts \big(\shape S {\mathcal{X}},\mathbf{Pyk}(\mathcal{S})\big) \longrightarrow {\Pyk(\mathcal{S})_{/K}}.
		\]
	induced by pre-composition with $ K \to \shape S {\mathcal{X}} $.
	By Lemma~\ref{Lemma: Precomp. pre. lim and colim}, the functor $ U $ admits a left adjoint $ L $.
\end{construction}

\begin{lemma}
	\label{prop:Lprecoh}
	Let $ P $ be a weakly contractible and compact object of $ \Pyk(\mathcal{S})_{/K} $.
	Then $ L(P) $ is a weakly contractible, compact and coherent object of  $\Fun \cts \big(\shape S {\mathcal{X}},\mathbf{Pyk}(\mathcal{S})\big)$.
\end{lemma}
\begin{beweis}
	Since the right adjoint $ U $ of $ L $ preserves colimits, it is clear that $ L(P) $ is weakly contractible and compact.
	Since $\mathcal{X}^{\pyk} \simeq \Fun \cts (\shape S {\mathcal{X}},\mathbf{Pyk}(\mathcal{S}))$ is locally coherent (see \cite[Theorem A.3.4.1]{lurie2016spectral}), it follows that we may find a collection $\{X_i\}_{i \in I}$ of coherent objects, such that the induced morphism
	\[
	\pi \colon \coprod_{i \in I} X_i \longrightarrow L(P)
	\]
	is an effective epimorphism.
	Since $L(P)$ is weakly contractible, there is a a section
	\[
	s \colon L(P) \longrightarrow \coprod_{i \in I} X_i
	\]
	of $\pi$ by Lemma~\ref{Lem:WContrSections}.
	Because $L(P)$ is furthermore compact, it follows that there is some finite $J \subseteq I$ such that~$s$ factors though $\coprod_{i \in J} X_i$ and thus $L(P)$ is a retract of the latter.
	Since finite coproducts of coherent objects are coherent, it follows that $L(P)$ is coherent by Lemma~\ref{Lemma: retracts of coh are coh}.
\end{beweis}

We are now ready to show the following:

\begin{satz}
	\label{thm:PyknotIsLocWeakCon}
	Let $\mathcal{X}$ be a spectral $\infty$-topos.
	Then $\mathcal{X}^{\pyk}$ is locally weakly contractible.
\end{satz}
\begin{beweis}
	By Theorem~\ref{Theorem: Cts Exodromy}, we equivalenty have to see that $ \Fun \cts (\shape S {\mathcal{X}},\mathbf{Pyk}(\mathcal{S}))$ is locally weakly contractible.
	By Lemma~\ref{prop:Lprecoh} it suffices to see that, for every $\mathcal{F} \in  \Fun \cts (\shape S {\mathcal{X}},\mathbf{Pyk}(\mathcal{S}))$, there is a collection of compact weakly contractible objects $P_i \in {\Pyk(\mathcal{S})_{/K}}$ and an effective epimorphism
	\[
		\coprod_i L(P_i) \longrightarrow \mathcal{F}.
	\]
	For this, we may pick a collection $P_i \in {\Pyk(\mathcal{S})_{/K}}$ and an effective epimorphism
	\[
		\beta \colon \coprod_i P_i\longrightarrow U(\mathcal{F}),
	\]
	which by adjunction corresponds to a morphism
	\[
		\alpha \colon \coprod_i L(P_i)  \simeq L\bigg( \coprod_i P_i\bigg)  \longrightarrow \mathcal{F}.
	\]
	We claim that $\alpha$ is an effective epimorphism.
	Because $U$ detects effective epimorphisms, it suffices to see that $U(\alpha)$ is an equivalence.
	But then we obtain a commutative triangle
	\[
	\begin{tikzcd}
\coprod_i P_i \arrow[rd, "\beta"'] \arrow[r, "\eta"] & U\bigg(L\bigg(\coprod_i P_i\bigg)\bigg) \arrow[d, "U(\alpha)"] \\
                                                      & U(\mathcal{F})                           
\end{tikzcd}		
	\]
	and since $\beta$ is an effective epimorphism, the map $U(\alpha)$ is an effective epimorphism as well, as desired.
\end{beweis}

\printbibliography

\end{document}